\documentclass[preprint]{elsarticle}
\usepackage{amsmath,array}
\usepackage{amssymb}
\usepackage{euscript}
\usepackage[all]{xy}
\numberwithin{equation}{section}
\usepackage[font=small,labelfont=bf]{caption}

\numberwithin{figure}{section}

\usepackage{amssymb,amsmath}
\usepackage{pstricks}
\usepackage{eepic}
\usepackage{times}

\newtheorem{definition}[equation]{Definition}
\newtheorem{theorem}[equation]{Theorem}
\newtheorem{proposition}[equation]{Proposition}
\newtheorem{corollary}[equation]{Corollary}
\newtheorem{lemma}[equation]{Lemma}

\newtheorem{example}[equation]{Example}

\newcommand{\ed}{{\rm d}}
\newcommand{\w}{{\mathchoice{\,{\scriptstyle\wedge}\,}{{\scriptstyle\wedge}}
      {{\scriptscriptstyle\wedge}}{{\scriptscriptstyle\wedge}}}}
\newcommand{\n}{\notag}
\newcommand{\noi}{\noindent}

\newcommand{\hook}{\hookrightarrow}
\newcommand{\tb}{\textbf}
\newcommand{\vsp}{\vspace}
\newcommand{\one}{\vspace{1mm}}
\newcommand{\two}{\vspace{2mm}}
\newcommand{\sq}{ \;  \square}

\newcommand{\tn}{\textnormal}

\newcommand{\C}{\mathbb{C}}
\newcommand{\R}{\mathbb{R}}
\newcommand{\I } {\mathcal{I}}

\newcommand{\W}{\mathcal{W}}
\newcommand{\A}{\mathcal{A}}
\newcommand{\PP}{\mathbb{P}}

\newcommand{\Co}{\mathcal{C}}

\newcommand{\F } {\mathcal{F}}
\newcommand{\x } {\textnormal{x}}

\newcommand{\Pf} {\noi\emph{Proof.}$\; \;$}

\newcommand{\beq}{\begin{equation}}
\newcommand{\eeq}{\end{equation} }

\begin{document}

\sloppy

\title{Maximum rank of a Legendrian web}
\tnotetext[tn1]{2000 Mathematics Subject Classification: 53A60}
\author{Joe S. Wang}
\ead{jswang12@gmail.com}
\address{KIAS, Seoul, South Korea}
%%%%%%%%%%%%%%%%%%%%%%%%%%%%%%%%%         address change
%\subjclass[2000]{53A60}
%\date{August, 2011}

\begin{abstract}
{We propose the Legendrian web in a contact three  manifold as a second order generalization of  the  planar web. An Abelian relation for a Legendrian web is analogously defined as an additive  equation among the first integrals of its foliations. For a class of  Legendrian $\, d$-webs defined by simple  second order ODE's, we give an algebraic construction  of
$$ \rho_d = \frac{(d-1)(d-2)(2d+3)}{6}$$
linearly independent Abelian relations. We then employ the method of local differential analysis and the theory of linear differential systems to show that $\, \rho_d$ is the maximum rank of a Legendrian $\, d$-web.

In the complex analytic category, we give a possible projective geometric interpretation of $\rho_d$ as an analogue of Castelnuovo bound for degree $2d$ surfaces in the 3-quadric  $\mathbb{Q}^3\subset\PP^4$ via the duality between $\PP^3$ and $\mathbb{Q}^3$ associated with the rank two complex simple Lie group $\tn{Sp}(2,\C)$.

The Legendrian 3-webs of maximum rank  three are analytically characterized, and their explicit local normal forms are found. For an application, we give an alternative characterization of a Darboux super-integrable metric as a two dimensional  Riemannian metric $\,g_+$ which admits a mate metric $\, g_-$ such that  a Legendrian
3-web naturally associated with the geodesic foliations of the pair  $\, g_{\pm}$ has maximum rank.
}
\end{abstract}

\maketitle
%-------------------------------------------------
%-------------------------------------------------
\section{Introduction}\label{sec1}
A  \emph{$\, d$-web  of  codimension $\,\nu$}  in a manifold is a set of  $\, d$  foliations  by  submanifolds   of  codimension $\, \nu$ whose $\, d$ tangent spaces are in general position at each point,  \cite{CG1}. Owing to the basic nature of its definition, webs  can be found in diverse areas of  differential  geometry,
\citep{BB,Chsurvey,CG2,GZ,PP,Aga,WaLeg}.

As  pointed out by Chern in the   survey \citep{Chsurvey}, Lie's theorem on the surfaces of double translation can be restated in web geometry terms as a result that a planar 4-web of  curves of   maximum  rank three (see below)   is algebraic, which  is also related to a  particular  case of the converse of Abel's theorem. For  an  example in projective  geometry, the classical  Chasles'  theorem   on the  intersecting  cubics is linked with the hexagonality of   algebraic planar 3-webs. According to \citep{Aga}, a three dimensional semi-simple  Frobenius manifold
carries a family of  characteristic hexagonal planar 3-webs. Planar 3-webs also  arise  as the asymptotic webs of  Legendrian surfaces in the projective space $\, \PP^5$, \citep{WaLeg}. Although in a possibly different context, the author's own interest in web geometry originated from Gelfand \&  Zakharevich's work on the Veronese webs attached  to  bi-Hamiltonian systems, \citep{GZ}. Much of  the developments of  web geometry up until 1930's are recorded in the book by Blaschke \&  Bol    ``Geometrie der Gewebe", \citep{BB}. For the more recent literature  on the subject, we refer to \citep{GS}\cite{PP} and the references therein.

\one
Let $\, \W$  be  a  $\, d$-web  of  codimension $\, \nu$. A single regular foliation (of fixed codimension) in a finite dimensional manifold has no local invariants, and the local geometry of   $\, \W$ lies entirely in the relative position among the set of $\, d$ foliations. A pertinent  notion in this aspect is the \emph{Abelian $\, k$-relation (equation)},  $\, 1 \leq k \leq  \nu$,
an additive relation  among the closed   basic $\, k$-forms with respect to each foliation.\footnotemark
\footnotetext{
An Abelian $\, 1$-relation represents an  additive relation among the first integrals of foliations.
}
The set of Abelian $\, k$-relations of  $\, \W$   naturally form  a vector space, and the \emph{$k$-rank} of $\, \W$ denoted by  $\, r^k(\W)$ is defined as the dimension of the vector space of  Abelian $\, k$-relations. The sequence of  integers $\,  r^k(\W)$, $\, k = 1, \, 2, \, ... \, \nu$, are the basic local invariants of a web.

\two
The $k$-rank of  a $\, d$-web   of codimension $\, \nu$ is not arbitrary. There generally exists an optimal bound
determined in terms of  $\,\,  k, \, d, \, \nu$, and the dimension of the ambient space, \citep{CG2}\citep{He}. It is a classical result of  Bol  that the maximum rank ($1$-rank) for a planar $\, d$-web is  $\, \frac{(d-1)(d-2)}{2}$, which is the well known genus bound for  a degree $\, d$ plane algebraic curve, \citep{Bol}. Chern   in his thesis   generalized this and showed that the maximum rank  for  a $\, d$-web of codimension 1
is equal  to the Castelnuovo bound for a degree $\, d$  algebraic curve,  \citep{Ch}.
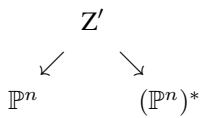
\begin{figure}
\begin{picture}(100, 50)(-80, 0)

\put(161,40){$\textnormal{Z}'$}
\put(175,25){$\searrow$}
\put(145,25){$\swarrow$}
\put(180,10){ $(\PP^n)^*$ }
\put(129,10){ $\, \PP^n$ }
\end{picture}

\caption{Projective duality}
\label{1duality}
\end{figure}

One of the main geometrical ideas behind  these results is the projective duality, Fig.\,\ref{1duality},
and the corresponding Abel's theorem, \citep{CG1}.\footnotemark
\footnotetext{
Many of  the arguments and results in web geometry
hold both in the real smooth, and the complex analytic categories.
Let us argue here in the complex analytic category.
}
A degree $\, d$ curve $\, \Sigma \subset (\PP^n)^*$ induces by the standard dual construction an associated $\, d$-web of  hyperplanes $\, \W_{\Sigma}$ in an open subset of $\, \PP^n$. By Abel's theorem,  the trace of a holomorphic  1-form on $\, \Sigma$ gives rise to an Abelian relation for $\, \W_{\Sigma}$. The dimension of the subspace of such algebraic   Abelian relations for $\, \W_{\Sigma}$ is then bounded by the Castelnuovo bound for the given algebraic curve.

With this background geometry, the results of  Bol, and Chern admit  the following  interpretation. The vector space  of Abelian relations  (for  a web of  codimension 1) can be analytically defined as the  space of solutions to a  linear differential system for the sections of a direct sum bundle canonically associated with a web. The differential system  is generally over-determined, see  Section \ref{sec2}. It is not difficult to imagine that such a uniform rank bound,  if exists, is determined  entirely by the symbol of the linear differential system. The results of  Bol, and Chern  indicate   that even if one considers the general webs of codimension 1, and not just the kind of algebraic webs  $\, \W_{\Sigma}$ described above, the relevant symbol of the linear differential system for Abelian relations is  isomorphic in an appropriate sense  to that of  the algebraic webs. The canonical linear differential system for Abelian relations does not get any more over-determined at the level of symbol when  one extends the category from algebraic webs to general webs.

\two
We  shall  consider in this paper  a generalization of web geometry to the geometric situation where the foliations are subject to a  first order constraint. Let   $\, X$ be a  three  dimensional  manifold. A  \emph{contact structure} on $\, X$ is a 2-plane field   $\, \mathcal{D} \subset TX$ locally defined by
$\, \mathcal{D} = \langle \, \theta \, \rangle^{\perp}$
for a 1-form $\, \theta$ that satisfies the non-degeneracy condition
$$ \ed \theta \w \theta \ne 0.$$
Given a contact structure, a \emph{Legendrian curve}  is an immersed $\, D$-horizontal curve. A \emph{Legendrian $\, d$-web} is  a set of $\, d$ pairwise transversal foliations by Legendrian curves. The Pfaff-Darboux  local  normal form theorem for contact structures  shows that a  Legendrian $\, d$-web can be locally described by a set of $\, d$ second  order ODE's.

Recall that a planar $\, d$-web is a set of $\, d$ transversal foliations by curves on a two dimensional manifold,
which is locally described by a set of $\, d$ first order ODE's. The purpose of this paper is to propose the Legendrian web  in a contact three manifold for  a second order generalization of  the  planar web. An  Abelian relation, which represents  an  additive functional relation  among the first integrals of   foliations, is  analogously defined for  the  Legendrian web. The rank of a Legendrian web is the dimension of the vector space of its Abelian relations.

\begin{figure}
\begin{picture}(100, 50)(-80, 0)

\put(161,40){$\textnormal{Z}$}
\put(175,25){$\searrow$}
\put(145,25){$\swarrow$}
\put(184,10){ $\mathbb{Q}^3$ }
\put(129,10){ $\, \PP^3$ }
\put(185,30){$\pi_2$}
\put(135,30){$\pi_1$}
\end{picture}

\caption{Projective duality  associated with  \,$\tn{Sp}(2,\C)$}
\label{1double}
\end{figure}
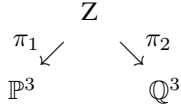

\one
In analogy with the discussion above,  in the complex analytic category the projective geometric model for our investigation is the  duality   associated with the rank two symplectic  group $\tn{Sp}(2,\C)$,  Figure\,\ref{1double}, see \citep{Br} for the details.
Here the 3-quadric
$\, \mathbb{Q}^3 \subset \PP^4$ is the space of  Legendrian lines in $\, \PP^3$,
and dually
$\, \mathbb{P}^3$ is the space of  null lines in $\,  \mathbb{Q}^3$
with respect to the $\tn{Sp}(2,\C)$ invariant  contact, and conformal structures   respectively.
The double fiber bundle $\, \textnormal{Z}$   is the  incidence space.

Consider a degree $\, 2d$  analytic surface  $\, \Sigma \subset \mathbb{Q}^3 \subset \PP^4$. Then generically  $\, \Sigma$ intersects a null line at $\, d$ distinct points, and it induces  by duality  an associated   $\, d$-web $\, \W_{\Sigma}$ of Legendrian lines on an open subset of $\, \PP^3$. This suggests to  consider  a  Legendrian  $\, d$-web
as a generalization of \emph{null degree $\, d$ surface in  $\, \mathbb{Q}^3 \subset \PP^4$}.\footnotemark
\footnotetext{Let us call a degree $\, 2d$ surface in  $\, \mathbb{Q}^3 \subset \PP^4$ simply a \emph{null degree $\, d$} surface.}

Let $\, \Omega$ be a closed meromorphic 1-form on a null degree $\, d$ surface $\, \Sigma$.
In consideration of Abel's theorem, one may declare that $\, \Omega$ is holomorphic
when the trace  $\,  (\pi_1)_* \circ  (\pi_2)^* \Omega$ vanishes.\footnotemark
\footnotetext{The relevant technical details are being ignored here. Since we are attempting to define   holomorphic objects on  generally singular surfaces, there are a few different choices, \citep{HP}. Note also that a holomorphic 1-form on a smooth compact complex surface is necessarily closed.}
By definition of trace, such a closed  holomorphic 1-form  gives rise to an Abelian  relation for $\, \W_{\Sigma}$. An Abelian relation for a Legendrian $\, d$-web can therefore be considered as a generalization of \emph{closed holomorphic 1-form on a null degree $\, d$ surface  in  $\, \mathbb{Q}^3$.}

\two
One of the initial motivations  for the present work came from the  observation that the bound for the  rank of  a Legendrian $\, d$-web would  imply the Castelnuovo  type bound for a null degree $\, d$ surface  in  $\, \mathbb{Q}^3$; determining  the   rank  bound by a direct local  analysis of  the linear differential system for Abelian relations, one may  derive  conversely a global Castelnuovo  bound for  the analytic surfaces in  $\, \mathbb{Q}^3$.

\two
\textbf{Main results.}

\one
1.
We give a local analytic characterization of  the Legendrian 3-webs of maximum rank three,\footnotemark
\footnotetext{Blaschke-Walberer Theorem states that
the maximum $1$-rank of  a  3-web of curves in a three dimensional manifold without any constraints is five,
and that  such 3-webs of maximum rank are essentially algebraic,
\cite[Section 35-37]{BB}.}
and determine their  local normal forms, Proposition \ref{propo56}, Theorem \ref{412thm}. As an application we give  an alternative characterization of a Darboux super-integrable metric in terms of a Legendrian 3-web associated with the geodesic flows of a pair of two dimensional Riemannian metrics, Theorem \ref{42thm}.

\one
2.
We give an algebraic construction of   $\frac{(d-1)(d-2)(2d+3)}{6}$ linearly independent Abelian relations
for a class of  simple   Legendrian $\, d$-webs, Section \ref{sec23}.

\one
3.
The  maximum rank  of a Legendrian $d$-web is, Section \ref{sec3},
\beq\label{mainresult}
\rho_d:= \frac{(d-1)(d-2)(2d+3)}{6}.
\eeq

\two
Let us mention some of the closely related works on web geometry.  In addition to the works of Bol, and Chern cited above, Chern \&   Griffiths   determined  the  bound for the $\, \nu$-rank of  the  webs  of  codimension $\, \nu$ in a manifold of dimension $\, \nu m$, \citep{CG2}. H\'enaut further generalized  Chern \&   Griffiths' bound to the bound for the $\, k$-rank for $\, 1  \leq  k \leq \nu$. We refer to \citep{Dam}\citep{Lit}\citep{Pir}\citep{Tr} for the analysis of the webs of  maximum rank.  Our main results can be considered as a  partial  generalization of   these works  to  the Legendrian web.  See also \cite{GL} for a recent survey on the rank problems  and linearizability for planar web.

Recently Shurygin gave an analysis of Legendrian 2 and 3 webs applying Kruglikov's description of differential invariants for the scalar second order ODE  up to point transformation, \cite{Shurygin1, Shurygin2}\cite{Kru}\cite{KL}.  

\two
In complex projective geometry the rank bound \eqref{mainresult}  calls for the study of  the extremal null degree $\, d$ surfaces in $\, \mathbb{Q}^3$ with the maximum number $\, \rho_d$ of (generalized) closed holomorphic 1-forms. The extremal null degree 3 surfaces are classified in \cite{Wa3}.

\two
We shall work within the real smooth category. Most of the results are valid in the complex analytic category with minor modification.  All of the arguments in this paper are local and they are concerned with the appropriate germs of geometric structures in a neighborhood of a reference point. The over-determined PDE machinery are used throughout the paper without specific references. For the standard references, we refer the reader to \citep{BCG3}\citep{IvLa}.

\two
We assume  the following well known property of  a  Vandermonde   matrix;
given a set of   $\, d$ distinct numbers $\,q^a, \, a=1, \, 2, \, ... \, d$,
the Vandermonde  matrix $\, ( \, V^b_a) =  (   \, (q^a)^{b-1})_{a,  b = 1}^d$ has full rank.

\two
We thank the anonymous referees for the valuable comments and suggestions which greatly improved the contents as well as the exposition of the paper. We thank Ben Mckay for the support and interest in this work.

%-------------------------------------------------
\subsection{Contact three manifold}\label{sec11}
This section recalls some basic definitions and properties of contact three manifolds.
The Pfaff-Darboux local normal form theorem (in dimension three) recorded in the below is a basic structure theorem which  gives a  local  uniformization of contact structures. For our purpose, it  provides an analytic representation   of  a  Legendrian web in terms of a  finite set of   second  order ODE's for one scalar function of one variable.

\two
Let $\, X$ be a three dimensional manifold. Let $\, TX \to X$ be the tangent bundle, and let  $\, T^*X \to X$ be the cotangent bundle of  $\, X$.
\begin{definition}\label{11defi}
A  \emph{\tb{contact  structure}} on a three dimensional manifold $\, X$ is a 2-plane field $\, \mathcal{D} \subset TX$
such that it is locally defined by
$\, \mathcal{D} = \langle \; \theta   \; \rangle^{\perp}$
for   a   1-form  $\, \theta$ which  satisfies   the non-degeneracy condition
\begin{align}\label{11nond}
\ed \theta   \w \theta   &\ne 0.
\end{align}
The \emph{contact line bundle} $\, \mathcal{C} \subset T^*X$
is the dual line bundle    $\, \mathcal{C} = \mathcal{D}^{\perp}$.
A  \emph{contact three  manifold} is a  three dimensional  manifold
with a  given  contact structure.
A \emph{(local)  contact transformation} of a contact three manifold
is a  (local)  diffeomorphism  which preserves the contact structure.
\end{definition}

The contact three manifold in the following example will be used throughout the paper.
\begin{example}\label{11example1}
The standard model of contact three manifold is the first jet space  $\, J^1(\R, \, \R)$ of scalar functions on the real line $\, \R$. Let $\, x$ be the standard coordinate of  $\, \R$, and let  $\, y, \, p\,$ denote the scalar function $\, y$ of $\,  x$, and the first  derivative   $\, p= \frac{\ed y}{\ed x}$ respectively. The set of three functions  $(x, \, y, \, p)$ form the standard coordinate of  $\, J^1(\R, \, \R)$. The contact line bundle is generated by the canonical contact 1-form
\beq\label{11canonical}
\theta = \ed y -p\, \ed x.
\eeq

Note the contact transformation $\, (x, \, y, \, p) \to (p, \, y-px, \,-x)$.
The full group of contact transformations is thus strictly larger than
the subgroup of contact  transformations induced by the point transformations of the form
$\, (x, \, y, \, p) \to
( \hat{x}(x, \, y), \,  \hat{y}(x, \, y),
\, \frac{\hat{y}_x+ p \hat{y}_y}{\hat{x}_x+ p \hat{x}_y})$.
\end{example}

\one
There exists a distinguished class of curves in a contact  three  manifold.
\begin{definition}
Let $\, X$ be a  contact  three   manifold
with the contact structure $\, \mathcal{D} \subset TX$.
A  \emph{\tb{Legendrian  curve}} in $X$ is an immersed $\, \mathcal{D}$-horizontal curve.
\end{definition}
\noi
From the non-degeneracy condition \eqref{11nond}, it is clear that there do not exist  immersed $\, \mathcal{D}$-horizontal surfaces.

\begin{example}[continued from Example \ref{11example1}]\label{11example2}
Let $f(x)$ be an arbitrary scalar function on $\R$.
By definition of the canonical contact 1-form \eqref{11canonical},
the first jet graph $$ \Big(x,\, f(x), \, \frac{\ed}{\ed x} f(x)=f'(x)\Big) \subset J^1(\R, \, \R)$$ is a Legendrian curve.
\end{example}

The following Pfaff-Darboux  local normal form theorem states that every Legendrian curve in a contact three manifold
is locally isomorphic to  Example \ref{11example2}.

\begin{theorem}[Pfaff-Darboux  local  normal form]\,
Let $\, X$ be a  contact  three   manifold.
Let $\, \x_0 \in X$ be a reference point.
There exists an adapted  local coordinate   $\, (x, y, p)$ of $\, X$
in a neighborhood of    $\, \x_0$  such that
the contact line bundle  is  generated  by
\beq\label{11normal}
\Co = \langle \, \ed y - p  \, \ed x \, \rangle.
\eeq
Let $\, \gamma \hook X$ be a  Legendrian curve through  $\, \x_0$.
Then  the adapted  local  coordinate  $\, (x, y, p)$ can be chosen so that
\beq\begin{array}{rll}
(x, y, p)&=(0,0,0)    & \mbox{ at}   \; \; \;  \x_0, \n \\
(y, p)&=(0,0)     &\mbox{\,along} \; \; \;  \gamma. \n
\end{array} \eeq
\end{theorem}
\Pf
The existence of  an  adapted local coordinate $\, (x, \, y, \, p)$ for \eqref{11normal} is well known, \cite[p38]{BCG3}.
Given the reference point $\, \x_0 \in X$,
let  $\, (x, \, y, \, p)|_{\x_0} =  (x_0, \, y_0, \, p_0)$.
From
\beq
\ed y - p \ed x  = \ed (y -p_0 x  - y_0 +p_0 x_0 ) - (p-p_0) \ed  (x-x_0), \n
\eeq
one may take
$\, ( \hat{x}, \,\hat{y}, \,\hat{p} )
= (x-x_0, \,y -p_0 x  - y_0 +p_0 x_0, \,p-p_0)$
as the desired coordinate  adapted at $\, \x_0$.

Let $\, (x, \, y, \, p)$ be such a coordinate system adapted at $\, \x_0$.
Given a Legendrian curve $\, \gamma$ through $\, \x_0$, there are three cases (they are not mutually exclusive):

$\bullet$\,
case $\, \ed x|_{\gamma} \ne 0$ at $\, \x_0$.

Let $\, \gamma^*y = f(x), \,  \gamma^*p = f'(x)$ locally
in a neighborhood of $\, \x_0$.
From
$\,
\ed y - p \ed x  = \ed (y - f(x) ) -  ( p -f'(x)) \ed x,
$
one may take
$\, ( \hat{x}, \,\hat{y}, \,\hat{p} )
= (x, \,y - f(x) , \,p-f'(x))$
as the desired coordinate  adapted   along $\, \gamma$.

$\bullet$\,
case $\, \ed y|_{\gamma} \ne 0$ at $\, \x_0$.

One may first translate $\, (x, \, y, \, p) \to (x, \, y+x, \,p+1)$
so that $\, (x, \, y, \, p)|_{\x_0} =  (0, \, 0, \, 1)$.
The contact 1-form $\, \ed y-p\ed x$
can then be written as
$\,
\ed y - p\ed x = -p (\ed x - \frac{1}{p} \ed y).
$
By interchanging $\, x$ and $\, y$, and up to a  translation similar as above,
this is reduced  to the first case.

$\bullet$\,
case $\, \ed p|_{\gamma} \ne 0$ at $\, \x_0$.

Under the contact transformation
$\, (x, \, y, \, p) \to (p, \, y-px, \,-x)$ from Example \ref{11example1}, this is also reduced to the first case.
$\sq$

\two
Pfaff-Darboux  local  normal form theorem  reflects the ampleness of  local group of contact transformations. We will see in the next section (Lie's theorem) that the local group of contact transformations is in fact large enough to normalize not just a single Legendrian curve but a foliation by Legendrian curves.

%-------------------------------------------------
%-------------------------------------------------
\section{Legendrian web}\label{sec2}
This section introduces the main subjects of the paper: Legendrian webs and their Abelian relations. In Section \ref{sec21}, a  Legendrian web  is defined as a finite set  of  foliations by Legendrian curves in a contact three manifold. In Section \ref{sec22}, an Abelian  relation for a Legendrian web is defined as an additive equation among the first integrals of the foliations.

%-------------------------------------------------
\subsection{Legendrian web}\label{sec21}
Let $d\geq 1$ be an integer.
\begin{definition}
A  \emph{\tb{Legendrian   $\, d$-web}}  on a  contact  three  manifold is a set of  $ \, d$  pairwise transversal  foliations by Legendrian curves.\footnotemark
\footnotetext{As remarked earlier, our analysis is essentially local and the issues such as global well-definedness of foliations shall be ignored.}
\end{definition}

Analytically a Legendrian web is described as follows.
Let $\, X$ be a  contact  three manifold with the contact line bundle $\, \Co  \subset T^*X$.
Let $\, \W$ be a Legendrian $\, d$-web on $\, X$. By duality, $\,\W$ is defined by a set of $\, d$ rank two sub-bundles
$\, \I^a \subset  T^*X, \, a = 1, \, 2, \, ... \, d$, such that $\,\Co\subset \I^a$ and from the transversality  condition that
$\, \I^a \cap \I^b = \Co$  for  $\, a \ne b$. This dual representation of a Legendrian web is convenient for the differential analysis for Abelian relations later on.

The Pfaff-Darboux local normal form theorem gives a concrete analytic formulation of a Legendrian web.
For a generic choice of adapted local coordinate $\, (x, \, y, \, p)$ for the given contact structure on $X$,
there exist a set of  $\,d$ functions
$\,   q^a(x,y,p), \, a=1, \, 2, \, ... \, d,$ $\, q^a \ne q^b$ for $\, a \ne b$,
such that the defining sub-bundle $\I^a$'s for $\W$ are locally generated by
$$
\I^a = \langle \, \ed y - p  \,\ed x, \,   \ed p - q^a \, \ed x \, \rangle,
\; \; a = 1, \, 2, \, ... \, d
$$
(and consequently on each Legendrian leaf we have $\ed x\ne 0$).
The geometry of  a Legendrian   $\, d$-web is thus locally equivalent to
the geometry of  a  set of  $\, d$ distinct second order ODE's
\beq\label{21ODE}
y^{''} =  q^a (x, \,y, \, y' ), \; a = 1, \, 2, \, ... \, d,\; q^a\ne q^b \;\, \mbox{for}\;\; a\ne b,
\eeq
up to contact transformation.

\two
Let us first examine the Legendrian $d$-webs for $d=1, 2$.

The following theorem of Lie states that there is no local invariant for a single Legendrian foliation and a Legendrian 1-web admits a unique local normal form up to contact transformation.
\begin{example}[Lie's theorem]
It is a classical theorem of Lie that a second  order ODE  for one scalar function of one variable
is locally  equivalent to the equation $\,  y^{''}=0$ up to contact transformation, \citep[p182]{Olv}.
The corresponding dual rank two sub-bundle $\, \I$ on $\,J^1(\R, \, \R)$ is
$$
\I  = \langle \, \ed y - p  \,\ed x, \,   \ed p   \, \rangle.$$
Under the contact transformation
$\, (x, \, y, \, p) \to (p, \, y-px, \,-x)$,
one has
\beq
\I  = \langle \, \ed y, \,   \ed x   \, \rangle. \n
\eeq
Lie's theorem thus shows that a Legendrian 1-web is locally equivalent to the foliation on $\,J^1(\R, \, \R)$ given by the fibers of canonical projection $\, J^1(\R, \, \R) \to J^0(\R, \, \R) \simeq \R^2$.
\end{example}

As a consequence this implies that the geometry of a Legendrian web lies in the relative position among the set of Legendrian foliations.

\two
A contact transformation of $\,J^1(\R, \, \R)$ which preserves the fibers of projection
$\, J^1(\R, \, \R) \to J^0(\R, \, \R) \simeq \R^2$
is locally induced from (the prolongation of) a point transformation of $\, \R^2$.
Lie's theorem above thus also implies that the geometry of  a Legendrian   $\, d$-web is locally equivalent to
the geometry of  a set of $\, d-1$ scalar second order ODE's
up to point transformation. In this regard, Cartan's work \citep{Ca} on the second order ODE's can be considered as a study of Legendrian 2-webs.\footnotemark
\footnotetext{See the argument below  Definition \ref{22Abeldefi}. A Legendrian 2-web has no nontrivial Abelian relations.}

\two
A class of Legendrian 3-webs are examined in detail in Section \ref{sec4}.

\two
Let us indicate in passing a geometrical construction for Legendrian webs. They are generalization to our setting of the well known construction in web geometry without constraints, \cite{PP}.

\begin{example}[continued from Example \ref{11example1}]\label{21family}
Let $\, J^1(\R, \, \R) \to J^0(\R, \, \R) \simeq \R^2$ be the canonical projection.
In terms of the standard adapted coordinate  $\, (x, \, y, \, p)$ of $\, J^1(\R, \, \R)$, the projection map is given by
$$ (x, \, y, \, p) \to (x, \, y).$$
Consider a two parameter family of  curves in $\, \R^2$,
each of which is locally a graph of a function $\, y = y(x)$.
Assuming this family of  curves  are   sufficiently generic to second order,
their  Legendrian   lifts $\, (x, \, y(x), \,y'(x) )$    to  $\, J^1(\R, \, \R)$
define  a Legendrian  foliation
on an open subset of  $\, J^1(\R, \, \R)$.
\end{example}
\begin{example}[continued]\label{21pencil}
Let $\, \{ \,  \textnormal{p}^1, \,  \textnormal{p}^2, \, ... \,  \textnormal{p}^d \, \} \subset \R^2$ be a set of $\, d$ distinct points. Consider for each $\,  \textnormal{p}^a$ the  two  parameter pencil of  $\, y$ as a quadratic  polynomial   of   $\, x$ based at $\,  \textnormal{p}^a$. A set of  $\, d$ distinct  points in $\, \R^2$ in this way gives rise to a Legendrian   $\, d$-web on an open subset of  $\, J^1(\R, \, \R)$.
\end{example}
\begin{example}[continued]\label{21LegBol}
Let $\, \{ \,  \textnormal{p}^1, \,  \textnormal{p}^2, \,  \textnormal{p}^3 \, \} \subset \R^2$ be a set of  three distinct  points. In addition to the quadratic pencil 3-web based at these points from Example \ref{21pencil}, consider the Legendrian foliation induced by  the two  parameter family of conics through the three points. A set of  three distinct  points  in $\, \R^2$ in this way gives  rise to a Legendrian  4-web on an open subset of  $\, J^1(\R, \, \R)$. This is  a Legendrian analogue of  Bol's exceptional 5-web, \citep{Dam}.
\end{example}

%-------------------------------------------------
\subsection{Abelian relation}\label{sec22}
In planar web geometry, the projective duality suggests  to consider a planar $\, d$-web as a  generalized degree $\, d$ plane curve. An Abelian relation for a planar web formally corresponds to a holomorphic 1-form through Abel's theorem, \citep{CG1}\cite{He}.

We shall follow this line of idea and give a definition of Abelian relation for a Legendrian web.
\begin{definition}\label{22Abeldefi}
Let $\, X$ be a contact three manifold. Let $\, \W$ be a Legendrian $\,d$-web on  $\, X$. Let $\, \I^a \subset  T^*X, \, a = 1, \, 2, \, ... \, d$, be the  dual rank  two  sub-bundles  which  define $\, \W$.
An \emph{\tb{Abelian relation}} of $\W$ is a $\, d$-tuple of 1-forms\footnotemark
\footnotetext{
In order to remove the ambiguity of  adding by constants, we  adopt here  the  differentiated  version of  first integrals for Legendrian foliation. See the remark below Proposition \ref{23propo} for an equivalent alternative definition.}
$\,(\Omega^1, \, \Omega^2, \, ... \, \Omega^d), \, \Omega^a \in  H^0(\I^a)$,
such that
\begin{align}
\ed\Omega^a    &= 0,    \; \; \, a = 1, \, 2, \, ... \, d,  \label{22defi} \\
\mbox{$\sum_b \Omega^b$} &= 0.  \label{22defi2}
\end{align}
The $\R$-vector space of  Abelian relations of $\W$  is denoted by $\A(\W)$.
The \emph{\tb{rank}} of a Legendrian web $\, \W$ is  the dimension of  $\,\A(\W)$.
\end{definition}

For a dual interpretation of Abelian relation for a Legendrian web as a generalized holomorphic 1-form, we refer to Figure\,\ref{1double} and the related remarks in Section \ref{sec1}. The idea to extend the notion of Abelian relation to Legendrian webs came from the exposition \cite{Br} on the dualities associated with the rank two simple Lie groups, see also Section \ref{sec5}.
Rank provides a basic local numerical invariant of a Legendrian web under contact transformation.

\two
Similarly as before, let us first examine the Abelian relations of Legendrian $d$-webs for $d=1, 2$.

A Legendrian 1-web  has rank zero by definition.

Let $\W$ be a Legendrian 2-web defined by a  pair of  rank two sub-bundles    $\, \I^a \subset  T^*X, \, a = 1, \, 2$. Let $\, (\Omega^1, \, \Omega^2)$ be an Abelian relation of $\W$. By definition we have $\, \Omega^1+\Omega^2=0$ and, since $\, \I^1 \cap \I^2 = \Co$, $\, \Omega^1$ must be a multiple of  contact 1-form. The defining  equation $\, \ed \Omega^1 = 0$ then forces $\, \Omega^1 = \Omega^2=0$.
Hence a Legendrian 2-web also has rank zero. Abelian relation is nontrivial for the Legendrian $\, d$-webs for $\, d \geq 3$. In view of the main theme of this paper, the  Abelian relations of a Legendrian web, we shall  restrict our attention to the Legendrian $\, d$-webs for the case 
$$\fbox{\;$d \geq 3$\;} $$  
from now on.

\two
The following example  shows that the Legendrian $d$-webs for $\,d\geq 3$ with at least one Abelian relation
exist in abundance.
\begin{example}\label{22atleastone}
Let $X$ be a contact three manifold with the contact line bundle $\Co \subset T^*X$. Let $\W_0$ be a   Legendrian $\, d$-web  ($d \geq 2$) defined by  the   $\, d$  rank  2 sub-bundles $\I^a \subset  T^*X, \, a = 1, \, 2, \, ... \, d$. Let $\, \Omega^a \in H^0(\I^a)$ be a closed 1-form for each $a = 1, \, 2, \, ... \, d$. Assuming  that $\, \Omega^a$'s are sufficiently generic, define  another distinct rank 2 sub-bundle $\, \I^{0} \subset T^*X$ by $\, \I^{0} = \Co  \oplus  \langle  \; \sum_b \Omega^b  \; \rangle$. The Legendrian $(d+1)$-web defined by $\I^a$'s,  $a = 0, \, 1, \, ... \, d$, possesses by construction at least one nontrivial Abelian relation.
\end{example}

One of the general problems in Legendrian web theory would be to give a characterization of Legendrian webs of maximal rank; the class of Legendrian webs with as many Abelian relations as possible.\footnotemark\footnotetext{The corresponding problem for planar $d$-webs remains open even for the first non-trivial case $d=5$, \cite{PP}\cite{Pir}.} Specifically, the defining equation for Abelian relations \eqref{22defi}, \eqref{22defi2} is a system of linear differential equations for the  sections of the direct sum bundle $\, \oplus_{a=1}^d \I^a$. A computation for the cases $\, d = 3, \, 4$ shows that it is over-determined, and one suspects that a generic Legendrian web has rank zero. The Legendrian webs we shall be interested in are those with maximal number of independent Abelian relations, and this imposes a stringent set of compatibility equations that such a Legendrian web must satisfy.  Moreover in this case it turns out that the vector space $\,\A(\W)$ inherits rich geometric structures from the layered feature of the integrability conditions for  \eqref{22defi}, \eqref{22defi2},
see Section \ref{sec3}. The geometry of $\,\A(\W)$ would be important in understanding the Legendrian webs of maximal rank.

On the other hand, it is not obvious if the rank of  a Legendrian $\,d$-web is bounded in terms of $\,d$, or even finite (for example, the characteristic variety  argument in \cite[p285]{BCG3} does not directly apply to  Legendrian  web). One of the basic problems would be to determine  an effective  bound on the rank, if such a bound exists. We shall show in Section \ref{sec3} that the rank of a Legendrian $\,d$-web admits the optimal bound

$$  \rho_d:= \frac{(d-1)(d-2)(2 d + 3)}{6}.$$

\vsp{1.5mm}

With regard to the duality associated with the rank two simple Lie group $\tn{Sp}(2,\C)$ mentioned earlier, this number should have a meaning in complex projective geometry as an analogue of Castelnuovo bound for surfaces in the 3-quadric $\mathbb{Q}^3\subset\PP^4$.

\two
Before we proceed to the differential analysis for the proof of rank bound, let us examine in some detail the first nontrivial case of Legendrian 3-webs. A direct computation gives an independent proof that the maximum rank of a Legendrian 3-web is  $\rho_3=3$. It also gives a taste of higher order differential analysis involved in analyzing the linear differential system for Abelian relations.

%-------------------------------------------------
%-------------------------------------------------
\section{Legendrian 3-webs of maximum rank}\label{sec4}
It is well known in planar web geometry that there exists essentially a unique local model for the 3-webs of maximum rank one up to diffeomorphism; a set of three families of parallel lines in the plane. As the analysis in this section will show, this uniqueness result is not true  for the case of Legendrian 3-webs and the local moduli space of Legendrian 3-webs of maximum rank three is generally two dimensional.\footnotemark
\footnotetext{A heuristic explanation for the dimension count of local moduli is as follows.
Consider the 3-quadric $\, \mathbb{Q}^3\subset \PP^4$.
The configuration space of three hyperplane sections in $\,\mathbb{Q}^3$ up to $\tn{Sp}(2,\C)$ motion has general dimension $\, 3\,\tn{dim}\,(\PP^4)^*  - \tn{dim}\,\tn{Sp}(2,\C)=3 \cdot 4 - 10 = 2.$ By projective duality, a  set of three hyperplane sections in $\mathbb{Q}^3$ corresponds to a linear Legendrian 3-web on an open subset of $\PP^3$, and this correspondence is generically injective. A direct computation shows that this class of Legendrian 3-webs have the maximum rank three. See \cite{Wa3} for the details.
}

In this section, we  employ the equivalence method of E. Cartan and  give a local analytic characterization of the Legendrian 3-webs of maximum rank,  Theorem \ref{412thm}. For an application,  we consider the Legendrian 3-web  associated  with a pair of two dimensional Riemannian metrics and give a new characterization of Darboux super-integrable metrics, Theorem \ref{42thm}.

%-------------------------------------------------
%-------------------------------------------------
\subsection{Local normal forms}\label{sec41}
Let us give a summary of the main results in this section.

A Legendrian 3-web $\W$ on a contact three manifold $X$ can be formulated as a $G_0$-structure on $X$ for a subgroup $G_0\subset\tn{GL}_3\R$; the collection of frames suitably adapted to the given 3-web $\W$ at each point of $X$ form a principal $G_0$ sub-bundle in the principal frame bundle of $X$, \eqref{410stgroup}. Under the normalization process of equivalence method this is further reduced to a $G\subset G_0$ structure with one dimensional structure group $G$. The torsion coefficients of the resulting structure equation and their successive derivatives are the basic local invariants of a Legendrian 3-web, Proposition \ref{propo56}.

The Abelian relations of $\W$ are by definition the solutions to a canonically attached linear differential equation. When expressed in terms of the adapted coframe of the reduced $G$-structure above, the compatibility conditions for $\W$ to have maximal rank are written as a set of linear equations in the basic local invariants of $\W$, \eqref{412compat}. A differential analysis of these equations shows that the torsion coefficients are necessarily covariant constant and that a Legendrian 3-web of maximum rank three is homogeneous for certain local Lie group structure, Theorem \ref{412thm}.

For readable references on the equivalence method, we refer to \cite{Ga}\cite{IvLa}.

%-------------------------------------------------
\subsubsection{Equivalence problem for Legendrian 3-webs}\label{sec411}
Let $\W$ be a Legendrian 3-web on a contact three manifold $X$.  Arguing locally, let $\theta$ be a contact 1-form, and let $\theta^a,\,a=1,2,3,$ be a set of 1-forms such that $\W$ is defined by the rank two sub-bundles $\,\I^a\subset T^*X$ locally generated by
\beq\label{41ideals}
\I^a=\langle\,\theta,\,\theta^a\,\rangle,\quad a=1,2,3.
\eeq
Up to scaling and adding multiples of  $\, \theta$, one may arrange so that
\beq\label{41defining}
\theta^1 + \theta^2 + \theta^3 =0.
\eeq

The most general transformation of $\,\{ \theta,\,\theta^a \}$ preserving the relation \eqref{41defining} is of the following form;
\beq\label{410stgroup}
 (\, \theta, \, \theta^a\,)  \to (\, \epsilon' \, \theta, \; \epsilon  \, \theta^{\tau(a)}+ t^a \, \theta\,),
\eeq
for  nonzero  scaling functions $\, \epsilon', \, \epsilon$, a permutation $\, \tau \in \mathcal{S}_3$ (the symmetric group on three letters), and translating factors $\, t^a$'s such that $t^1 + t^2 + t^3 = 0$.
This can be utilized to obtain a refined set of generators.
\begin{lemma}\label{41lem}
There exists a transformation of 1-forms $\,\{ \theta,\,\theta^a \}$ as described above such that
\begin{align}\label{41refined}
  \;\;\,\;\; \ed\theta &\equiv \theta^1\w\theta^2\mod\theta, \\
  \;\;\; \ed\theta^a&\equiv 0\qquad \,\mod \theta,\quad a=1,2,3.\n
\end{align}
\end{lemma}
\Pf
For a given $\,\{ \theta,\,\theta^a \}$ let
\begin{align}
\ed\theta &\equiv s\theta^1\w\theta^2\;\; \mod\theta,\n\\
\ed\theta^a&\equiv s^a\theta^1\w\theta^2 \mod \theta,\quad a=1,2,3,\n
\end{align}
for coefficients $s, \,s^a$. Since $\theta$ is a contact 1-form, $s$ is a nonzero function. Consider the transformation
\begin{align}
\theta &\to s\theta,\n \\
\theta^a&\to \theta^a-\frac{s^a}{s}\theta. \n \qquad\qquad\sq
\end{align}

Note that for the refined set of 1-forms $\,\{ \theta,\,\theta^a \}$ in Lemma \ref{41lem},
the most general transformation preserving the relation \eqref{41refined} is 
\beq\label{411stgroup}
 (\, \theta, \, \theta^a\,)  \to (\, \mbox{sgn}(\tau) \epsilon^2 \, \theta, \; \epsilon  \, \theta^{\tau (a)})
\eeq
for a nonzero  scaling function $\, \epsilon$,
and a permutation $\, \tau \in \mathcal{S}_3$.

\two
The analysis above can be recast in the language of Cartan's equivalence method as follows.
Let $\W$ be a Legendrian 3-web on a contact three manifold $X$. Let $F \to X$ be the principal $\tn{GL}_3\R$-frame bundle. With an abuse of notation, let $\{\theta, \theta^1, \theta^2\}$ denote the tautological 1-forms on $F$ and set
$\theta^3=-(\theta^1+\theta^2)$.
Then, there exists a principal sub-bundle $B\subset F$ with one dimensional structure group $G$ acting on $\{\theta, \theta^1, \theta^2\}$ by \eqref{411stgroup} such that;
for any local section $\sigma:X\to B$

\one
\qquad a) $\sigma^*\theta$ is a contact 1-form, and $\W$ is defined by the ideals generated by $\I^a=\{\sigma^*\theta, \sigma^*\theta^a\}$,  $a=1,2,3$,

\qquad b) the 1-forms $\{\sigma^*\theta, \sigma^*\theta^a\}$ satisfy the relations \eqref{41defining}, \eqref{41refined}.

\one
From the general theory, it is not difficult to show that the sub-bundle $B$ can be chosen so that the tautological 1-forms  satisfy the following structure equation on $B$ for a choice of pseudo-connection 1-form $\alpha$. We leave the details to the reader.

\beq\label{41struct}
\ed
\begin{pmatrix}
\theta^1 \\ \theta^2 \\  \theta
\end{pmatrix}
= -
\begin{pmatrix}
\alpha & \cdot &  \cdot \\
\cdot &  \alpha & \cdot \\
\cdot &  \cdot &  2 \alpha
\end{pmatrix}
\w
\begin{pmatrix}
\theta^1 \\ \theta^2 \\  \theta
\end{pmatrix}
+
\begin{pmatrix}
\theta \w ( R \, \theta^1 + S \, \theta^2 )  \\
\theta \w ( T \, \theta^1 - R\, \theta^2 )  \\
\theta^1 \w \theta^2
\end{pmatrix}.
\eeq
Here '$\cdot$' denotes 0, and $\, R, \, S, \, T$ are torsion coefficients.\footnotemark
\footnotetext{
The structure group of the principal bundle $ B$ acts on
$(\theta,   \theta^a)$  by  
$ ( \theta,  \theta^a\,) \to(  \mbox{sgn}(\tau) \epsilon^2 \, \theta, \, \epsilon \, \theta^{\tau (a)} \,)$
for  a nonzero  scaling parameter $\, \epsilon$, and a permutation by $\, \tau \in S_3$. For example, under the permutation
$\, (\,\theta^1, \, \theta^2, \, \theta^3) \to
    (\,\theta^2, \, \theta^1, \, \theta^3)$,
the torsion coefficients transform to
\beq\label{41permute}
\begin{pmatrix}
R & S \\
T & -R
\end{pmatrix}
\to
\begin{pmatrix}
R &  -T \\
 -S &  -R
\end{pmatrix}.
\eeq}

Note that $\alpha$ is uniquely determined by this structure equation.
The exterior derivative identity $\, \ed(\ed(\theta))=0$ furthermore implies that
\beq\label{41alpha}
\ed \alpha = \theta \w ( \,N \theta^1 + L \, \theta^2 \, ),
\eeq
for coefficients $\, N, \, L$.

\two
Let us summarize the analysis so far.
\begin{proposition}\label{propo56}
Let $\, \W$ be a Legendrian 3-web on a contact three manifold $\, X$. 
There exists an adapted sub-bundle $\, B$ of the $\,  \textnormal{GL}_3 \R$ principal frame bundle  of $\, X$ on which the tautological 1-forms $\,  \theta, \, \theta^a , \, a = 1, \, 2, \, 3$, satisfy the equations \eqref{41ideals}, \eqref{41defining}, \eqref{41struct}, \eqref{41alpha}. The functional relations among the torsion coefficients $\, R, \, S, \, T, \, N, \, L$ and their successive  derivatives are the basic local invariants of $\, \W$.

The algebra of  local symmetry vector fields of $\, \W$ is at most four dimensional, and it is four dimensional whenever  the torsion coefficients $\, R, \, S, \, T, \, N, \, L$ all vanish.
\end{proposition}
\Pf
The  algebra of  local symmetry vector fields of   $\, \W$  is  four dimensional
when the torsion coefficients $\, R, \, S, \, T, \, N, \, L$ are all constants.
Since they scale along the fibers of the bundle $\, B$,
e.g., $\, d R \equiv  2 R \alpha, \mod \; \theta, \, \theta^1, \, \theta^2$,
this occurs only when $\, R, \, S, \, T, \, N, \, L$ all vanish.
The rest  follows from the general theory of equivalence method.
$\sq$

\one
The differential analysis for Abelian relations in the next sub-section will be carried out on the bundle $B$.

\two
We record in passing the following corollary of Proposition \ref{propo56} on the intersection of local  algebra of  symmetry vector fields of two distinct path geometries on a surface.

Let $\, M$ be a two dimensional manifold. Let $\, \PP(TM) \to M$ be the projective tangent bundle
equipped with the canonical contact structure.
Recall  that  a path geometry on $\, M$ is a (local) Legendrian foliation on $\, \PP(TM)$
transversal to the fibers of projection  $\, \PP(TM) \to M$,  \citep{Brpath}.
\begin{corollary}
Let $\, \mathcal{F}_{\pm}$ be a pair of  distinct  path geometries on a surface  $\, M$.
Let $\, \mathfrak{p}_{\pm} \subset   \mathfrak{X}(M)$ be the algebra of  symmetry vector fields for  $\, \mathcal{F}_{\pm}$ respectively (here $\,\mathfrak{X}(M)$ is the algebra of vector fields on $\, M$). Then
\beq
\textnormal{dim}\;   \mathfrak{p}_+ \cap \mathfrak{p}_- \leq   4. \n
\eeq
The equality holds whenever  the pair $\, \mathcal{F}_{\pm}$ is locally equivalent to the pair of flat path geometries
defined by the second order ODE's
\beq
y^{''} = q_{\pm}, \quad \; \;     q_+,    \, q_- \; \textnormal{are  distinct constants}. \n
\eeq
\end{corollary}
\Pf
The  prolongation of  the intersection  $\;  \mathfrak{p}_+ \cap \mathfrak{p}_-$ to   $\, \PP(TM)$
is the   symmetry algebra of vector fields for the Legendrian 3-web defined by $\, \mathcal{F}_{\pm}$, and the fibers of projection  $\, \PP(TM) \to M$.
$\sq$

%-------------------------------------------------
%-------------------------------------------------
\subsubsection{Abelian relations}\label{sec412}
We continue the analysis from Proposition \ref{propo56}.

\two
For three functions $f, g^1, g^2$ on $B$, set
\begin{align}
\Omega^1 &= f \, \theta^1 + g^1 \, \theta, \n \\
\Omega^2 &= f \, \theta^2 + g^2 \, \theta, \n \\
\Omega^3 &= - \Omega^1 - \Omega^2, \n
\end{align}
denote a section of $\, \I^a$ for $a = 1, \, 2, \, 3$, which satisfy  the relation $\, \sum_b \Omega^b=0$.
The covariant derivatives of $\, f, \, g^1, \, g^2$ are written by
\begin{align}
\ed f &= f\, \alpha + f_1 \theta^1 + f_2 \theta^2 + f_0 \theta, \n \\
\ed g^1 &= 2 g^1\,   \alpha + g^1_1 \theta^1 + g^1_2 \theta^2 + g^1_0 \theta, \n \\
\ed g^2 &= 2 g^2\,   \alpha + g^2_1 \theta^1 + g^2_2 \theta^2 + g^2_0 \theta. \n
\end{align}
We shall adopt the  similar notational convention for the covariant derivatives of the structure coefficients
$\, R, \, S, \, T, \, N, \, L$, e.g.,
$\, \ed R = 2 R \, \alpha + R_1 \, \theta^1 + R_2 \, \theta^2 + R_0 \, \theta$, etc.

By direct computation one finds that the defining equation for Abelian relation 
$$\,\ed  \Omega^a = 0, \quad a = 1, \, 2,$$
is equivalent to that
\beq\label{42eq1}\begin{array}{ll}
f_2 =g^1,  &f_1 = - g^2,   \\
g^1_2 =S f, &g^2_1 = T f,   \\
g^1_1=f_0 + R f, &g^2_2 = f_0 - R f.
\end{array}
\eeq
With these  relations imposed,
the equation $\,\ed (\ed (f)) \w \theta = 0$ implies $\, f_0 =0.$
Differentiating again, the exterior derivative identity $\, \ed (\ed (f))  = 0$ gives
\begin{align}\label{42eq2}
g^1_0 &= - L f+ S g^2 +R g^1,   \\
g^2_0 &=\; \, N f -R g^2 + T g^1. \n
\end{align}
The linear differential system for the coefficients
$\, \{ \, f, \, g^1, \, g^2 \, \}$
closes up at this step. It follows from the theory of linear differential equations that a Legendrian 3-web admits at most three linearly independent Abelian relations.

\two
We proceed to examine the consequences of \eqref{42eq1}, \eqref{42eq2}.
A direct computation shows that the identities  $\, \ed (\ed (g^1)) = 0, \, \ed (\ed (g^2)) = 0$ give the following three compatibility equations.
\begin{align}\label{412compat}
(S_0+L_2) f +(-S_1+ 4 L ) g^1 - S_2 g^2 &=0,   \\
(R_0 + L_1) f +(2 N- R_1) g^1 +(-S_1-L) g^2 &=0, \n \\
(-T_0 +N_1 ) f + T_1 g^1 +(-3 N-R_1 ) g^2 &=0. \n
\end{align}
In order for a Legendrian 3-web to have the maximal rank three,
these equations should hold identically and each of the nine coefficients of $\, f, \, g^1, \, g^2$ in \eqref{412compat}
must vanish,
otherwise the rank is strictly less than three.
This implies, after a short computation,
that the structure coefficients of a Legendrian 3-web must satisfy
\begin{align}\label{42const}
\ed R&\equiv \ed S \equiv \ed T \equiv 0,   \mod \; \alpha,   \\
N&=L=0. \n
\end{align}
Conversely, it is clear that when these equations hold a Legendrian 3-web admits three dimensional space of Abelian relations.

\one
The question remains regarding the existence of Legendrian 3-webs satisfying \eqref{42const}.
The structure equation \eqref{41struct}, \eqref{41alpha} with \eqref{42const} imposed is easily seen to be compatible, i.e., $\ed^2=0$ is a formal identity of the structure equation, and it follows from the general theory of $G$-structures, \cite[Appendix]{BrBoch}, that there exist two parameter family of Legendrian 3-webs with the structure coefficients satisfying \eqref{42const}. In particular, \emph{the maximum rank of a Legendrian 3-web is $\rho_3=3$.}

\two
We wish to integrate the structure equation \eqref{41struct}, \eqref{41alpha} with \eqref{42const} and determine  a local normal form for the Legendrian 3-webs of maximum rank in terms of a pair of second order ODE's. Since $\ed\alpha=0$,  for simplicity let us take a section of the frame bundle $B\to X$ for which $\, \alpha= 0$.
Then $\, R, \, S, \, T$ become constants (up to nonzero scale depending on the choice of section), and \eqref{41struct} is the  structure equation of a Lie group.\footnotemark
\footnotetext{See \cite[p29]{BrLie} for a classification of
three dimensional Lie groups.}\,
Hence a Legendrian 3-web of maximum rank three is necessarily locally homogeneous.
\begin{theorem}\label{412thm}
Let $\, X$ be a contact three manifold. Let $\, \W$ be a Legendrian 3-web on $\, X$. 

\one a)
The maximum rank of $\W$ is three.

\one b)
Suppose $\W$ has the maximum rank. Then the structure invariants of $\, \W$ satisfy the relation \eqref{42const}, and $\, \W$ is locally equivalent to a left invariant homogeneous Legendrian 3-web on a three dimensional Lie group.

\one cexpressions below)
There exists a local contact isomorphism of $ X$ with the 1-jet space $ J^1(\R, \, \R)$ equipped with the canonical contact structure such that $\, \W$ is defined by the integral curves of  the following set of second order ODE's. Note that $\, R, \, S, \, T$'s in the expressions below are the (constant) structure coefficients of  \eqref{41struct} with $\, \alpha = 0$.

%$\bullet$
%\texttt{Case $\,R^2 + S T=0, \, R \ne 0$.}\;
%\begin{align}
%y^{''}+  T y &=0, \n \\
%y^{''}+  T y - \frac{T}{R} &=0, \n \\
%y^{''}+  T y - \frac{T^2}{R (R+T)} &=0. \n
%\end{align}
%When $\, R+T=0$,
%the Legendrian foliation by the third ODE
%degenerates to the fibers of the natural projection
%$\, J^1(\C, \, \C) \to J^0(\R, \, \R)$.
\one
$\bullet$
\textbf{Case $\, R^2 + S T=0$.}\;
When $\, R =S = 0$,
\begin{align}
y^{''}+ T  y &=0, \n \\
y^{''}+ T  y + 1 &=0, \n
\end{align}
and the Legendrian foliation by the fibers of  projection $\, J^1(\R, \, \R) \to J^0(\R, \, \R)$.
The general cases are obtained by applying
the linear action \eqref{42conjugation}
to \eqref{42T}.

$\bullet$
\textbf{Case $\,R^2+ST > 0$.}\;
When $\, R >0, \, S = T = 0$,
\begin{align}\label{ODERplus}
y^{''}+ R (y')^2 &=0, \n  \\
y^{''}- R (y')^2 &=0, \n \\
y^{''}+ R \tanh{(R y)} (y')^2 &=0. \n
\end{align}
The general cases are obtained by applying
the linear action \eqref{42conjugation}
to \eqref{42R}.

$\bullet$
\textbf{Case $\,R^2+ST < 0$.}\;
When $\, R = 0, \, -S= T \ne  0$,
%When $\, R = 0, \, -S= T < 0$,
%\begin{align}
%y^{''}  -  T \tan(Ty)  (y')^2 &=0, \n \\
%y^{''}  + T \cot(Ty)  (y')^2 &=0, \n \\
%y^{''}  + T \frac{   \cos(Ty)-\sin(Ty)  }{ \cos(Ty) + \sin(Ty)}  (y')^2 &=0. \n
%\end{align}
%When $\, R = 0, \, -S= T \ne  0$,
\begin{align}
y^{''}  +   \tan ( T y ) ( 1 - T  ( y' ) ^{2}) &=0, \n \\
y^{''}  -    \cot ( T y ) ( 1 - T  ( y' ) ^{2}) &=0, \n \\
y^{''}  -\frac{   \cos(Ty)-\sin(Ty)  }{ \cos(Ty) + \sin(Ty)}( 1 - T  ( y' ) ^{2})&=0. \n
\end{align}
The general cases are obtained by applying
the linear action \eqref{42conjugation}
to \eqref{42mTplus}.
\end{theorem}

\noi
\Pf
The proof is by direct computation.
Set the contact 1-form $\, \theta = dy-p dx$.
Then $\, \theta$ and $\, \theta^1, \, \theta^2$ below
satisfy  the structure equation \eqref{41struct},  \eqref{42const}
with $\, \alpha = 0$.
Note that under the linear action by an element $\, \textnormal{g} \in \textnormal{SL}_2\R$,
\beq\label{42conjugation}
\begin{pmatrix}
\theta^1 \\  \theta^2
\end{pmatrix}
 \to \textnormal{g}
\begin{pmatrix}
\theta^1 \\  \theta^2
\end{pmatrix}, \quad
\textnormal{then} \quad
\begin{pmatrix}
R & S \\
T & -R
\end{pmatrix}
 \to  \textnormal{g}
\begin{pmatrix}
R & S \\
T & -R
\end{pmatrix}
\textnormal{g}^{-1}.
\eeq

$\bullet$
\textbf{Case $\,R^2 + ST=0$}.\;
Up to the conjugation  \eqref{42conjugation},
one may assume
$\, R = S = 0.$
Set
\begin{align}\label{42T}
\theta^1&=   \ed x,  \\
\theta^2&=   \ed p + T  y \ed x, \n \\
\theta^1+\theta^2&=
  \ed p + (T  y + 1) \ed x. \n
\end{align}

$\bullet$
\textbf{Case $\,R^2+ST > 0$}.\;
Up to  \eqref{42conjugation},
one may assume $\, R > 0, \,  S = T = 0$.
Set
\begin{align}\label{42R}
\theta^1&= \frac{\exp{(R y)}} {\sqrt{2R} \, p} (\ed p + R p^2  \ed x),   \\
\theta^2&= \frac{\exp{(-R y)}}{\sqrt{2R} \, p} (\ed p - R p^2  \ed x), \n \\
\theta^1+\theta^2&=
\frac{2 \cosh{(R y)}}{\sqrt{2R} p} (\ed p + R \tanh{(R y)} p^2  \ed x).  \n
\end{align}

$\bullet$
\textbf{Case $\,R^2+ST < 0$}.\;
Up to   \eqref{42conjugation},
one may assume $\, R=0.$
%In case $\, T < 0$,
%set
%\begin{align}\label{42mTminus}
%\theta^1&= \frac {\cos ( Ty ) }{\sqrt {-T}p}
%                   (dp - T  \tan(T y) p^2 dx),      \\
%\theta^2&= \frac {\sin ( Ty )  }{\sqrt {-T}p}
%                   (dp + T  \cot(T y)  p^2 dx),  \n   \\
%\theta^1+\theta^2&=
%\frac {\cos ( Ty ) +\sin  ( Ty ) }{\sqrt {-T}p}
%(dp   +  T  \frac { \cos ( Ty ) - \sin  ( Ty )  }{\cos ( Ty ) +\sin  ( Ty )} p^2 \, dx).  \n
%\end{align}
%%%%%%%%%%%%%%%%%%%%%%%%%%%%%
%In the general case,
Set %(this  also works for the case $\, T < 0$)
\begin{align}\label{42mTplus}
\theta^1&={\frac {\cos  ( Ty ) }{\sqrt {1-T{p}^{2}}}}
 (
\ed p + \tan(Ty) (1 - T p^2    ) \ed x
 ),    \\
\theta^2&={\frac {\sin  ( Ty  ) }{\sqrt {1-T{p}^{2}}}}
 (
\ed p - \cot(Ty) (1 - T p^2    ) \ed x
 ),    \n    \\
\theta^1+\theta^2&=
{\frac {\cos ( Ty )  + \sin  ( Ty  ) }{\sqrt {1-T{p}^{2}}}}
 (
\ed p - \frac{ \cos ( Ty )  -  \sin  ( Ty  ) }{\cos ( Ty )  + \sin  ( Ty  ) } (1 - T p^2    ) \ed x
 ).    \n
\sq
\end{align}
%-------------------------------------------------
%-------------------------------------------------
\subsection{Geodesic Legendrian webs of maximum rank}\label{sec42}
There is a geometric situation where Legendrian webs naturally occur. Let $M$ be a two dimensional surface.  Let $\PP(TM) \to  M$ be the projective tangent bundle equipped with the canonical contact structure. Consider on $M$ a finite set of Riemannian metrics $\{ \,g^a\,\}_{a=1}^d$. Each metric $g^a$ defines the geodesic Legendrian foliation on $\, \PP(TM)$. Combined with the foliation by fibers of projection $\PP(TM) \to  M$, which are Legendrian, a set of $d$  Riemannian metrics on a surface $M$ gives rise to a Legendrian $(d+1)$-web $\W_{\{g^a\}}$ on $\PP(TM)$.

In this section, we give  an application of the preceding analysis to the case a Legendrian 3-web is   defined  by  a pair  of  two dimensional Riemannian metrics.
\begin{definition}\label{42defi}
Let $\,  g_{\pm}$ be a pair of Riemannian metrics on a two dimensional surface $\, M$.
Under an appropriate transversality condition,\footnotemark
\footnotetext{See \textbf{Step 1} in Section \ref{sec421}.
}
let $\, \W_{g_{\pm}}$ be the Legendrian 3-web on $\, \PP(TM)$ defined by the geodesic foliations of $\,  g_{\pm}$ and the fibers of projection  $\, \PP(TM) \to  M$. The pair of Riemannian metrics are \emph{\tb{maximally geodesically compatible}} when the associated Legendrian 3-web $\, \W_{g_{\pm}}$ has the maximum rank three. In this case, $\,  g_{\pm}$ are the \emph{maximally geodesically compatible mates}
to each other.
\end{definition}

Theorem \ref{412thm} implies that  a maximally geodesically compatible pair of metrics share at least three dimensional algebra of projective vector fields. The classification result of  \citep{Bret} then immediately gives the local normal forms for such pair of metrics, Theorem \ref{42thm}.

\two
Before we state the main result, let us  recall the relevant known facts on the two dimensional Riemannian metrics.
We refer to \citep{Bret}\citep{Mat}\citep{Mat2} for the details.

Let $\, M$ be a two dimensional oriented Riemannian manifold equipped with the metric $\, g$. Let $\, \pi: \, FM \to M$ be the $\, \textnormal{SO}_2$-bundle of oriented orthonormal frames. Let  $\, \{ \, \omega^1, \, \omega^2 \, \}$ be the tautological 1-forms on $\, FM$ such that  $\, \pi^*g = (\omega^1)^2 + (\omega^2)^2$. The Levi-Civita connection 1-form $\, \rho$, and the Gau\ss\, curvature $\, K$  are uniquely defined on $FM$ by the equations
\begin{align}\label{42struct}
\ed \omega^1 &= -  \rho \w \omega^2,   \\
\ed \omega^2 &= \; \; \; \rho \w \omega^1, \n \\
\ed \rho          &= K \, \omega^1 \w \omega^2. \n
\end{align}
Differentiating the last equation, the higher  order derivatives of $\, K$ are inductively defined by
\begin{align}
\ed K &= K_1 \omega^1 + K_2 \omega^2, \n \\
\vspace{5mm}
\ed K_1 &= -K_2 \rho      + K_{11} \omega^1 + K_{12} \omega^2, \n \\
\ed K_2 &= \; \; \; K_1 \rho + K_{21} \omega^1 + K_{22} \omega^2,   \quad K_{12} = K_{21},
\quad  \textnormal{etc}. \n
\end{align}

\two
The set of unparameterized, oriented geodesics of $\, g$
defines the geodesic flow (foliation)  $\, \mathcal{F}_g$ on $\, \PP(TM)$,
which is  naturally extended to the tangent bundle  $\, TM  \to M$.

Let $\, \odot^k(T^*M) \to M$ be the bundle of symmetric $\, k$-forms on $\, M$,
which we consider as the functions on $\, TM$ homogeneous of degree $\, k$ on the fibers.
Set
\begin{align}
I_k(g)  &= \, \{ \, \sigma \in H^0(\odot^k(T^*M)) \;\; | \; \; \sigma \; \;
                     \textnormal{is a first integral for} \;  \mathcal{F}_g  \; \; \}, \n \\
\varrho_k(g) &=\textnormal{dimension}_{\R}  \,   I_k(g) . \n
\end{align}

The following results on  $\, I_k(g)$ for $\, k = 1, \, 2, \, 3$,   are known.
The class of  Riemannian  metrics we shall be interested in are  the Darboux super-integrable metrics
for which $\, \varrho_2(g) =4$.

\vsp{2mm}

$I_1(g)$: \\
An element  of  $\, I_1(g)$ corresponds to a  Killing vector field  of $\, g$.
For any metric $\, g$, one necessarily has
$\, \varrho_1(g) = 3, \,   1$, or $\, 0$ (generic case);

$\, \varrho_1(g) =3$ when $\, g$ is a metric of constant curvature.

$\, \varrho_1(g) =1$ when $\, g$ admits a Killing vector field, which is unique up to scale.
For an analytic characterization of this class of metrics, \cite[p321]{Cartan}.
For instance, the curvature must satisfy the equation
 $\, K_1 K_2 (K_{11} - K_{22}) - (K_1^2-K_2^2) K_{21} = 0$.

\vsp{2mm}

$I_2(g)$: \\
An element  of  $\, I_2(g)$ corresponds to a Riemannian metric projectively equivalent to $\, g$,
\citep{Mat2}.
For any metric $\, g$,
$\, \varrho_2(g) = 6, \, 4, \, ... \,$, or $\, 1$(generic case);

$\, \varrho_2(g) =6$ when $\, g$ is a metric of constant curvature.

$\, \varrho_2(g) =4$ when $\, g$ admits a Killing field  ($\, \varrho_1(g) =1$),
and the curvature of $\, g$ satisfies an additional  set of second, and third order differential equations.
This class of metrics are called \emph{Darboux super-integrable}, \citep{Bret}.

Darboux super-integrable metrics are equivalently characterized by
having three   dimensional algebra of  projective vector fields.
Combining this with a classical result of Lie
on the classification of the algebra of projective vector fields on the plane,
\citep{Bret} determined the explicit local normal forms
for this class of metrics.

\vsp{2mm}

$I_3(g)$: \\
The cubic integrals are studied in \citep{Mat}.
For any metric $\, g$,
$\, \varrho_3(g)  \leq 10$,
and
$\, \varrho_3(g) =10$ when $\, g$ is a metric of constant curvature.
It is possible  that the next admissible value of $\,\varrho_3(g) $ is $\, 4$,
which is attained  by  the  Darboux super-integrable metrics.

\two
Let us now state the main results of this section.
\begin{theorem}\label{42thm}
Let $\,  g_{+}$ be a  two dimensional Riemannian metric. $\, g_+$ admits a maximally geodesically compatible mate $\, g_-$ whenever  $\, g_+$ is  either of constant curvature, or Darboux super-integrable. 

\one
Let $\, \mathcal{R}_{g_+}$ denote the moduli space of  Riemannian metrics maximally geodesically compatible with the given metric $\, g_+$. 

\;a)
when $\, g_+$ is of constant curvature,
$\, \mathcal{R}_{g_+}$ consists of  10 parameter family of metrics of constant curvature,
and
8  parameter family of  Darboux super-integrable metrics.

\;b)
when $\, g_+$ is Darboux super-integrable,
$\, \mathcal{R}_{g_+}$ consists of  5 parameter family of metrics of constant curvature,
and
4 parameter family of  Darboux super-integrable metrics.
\end{theorem}

Many of the metrics in this  theorem  are geodesically equivalent.
For example,
for the class of Darboux super-integrable metrics $\, g_+$
to be considered in Section \ref{sec422},
all of  the  5 parameter family of  mates  of constant curvature  are geodesically equivalent,
and the 4 parameter family of  Darboux super-integrable mates are foliated by
        3 parameter families of  geodesically equivalent metrics.

\subsubsection{Differential analysis}\label{sec421}
The differential analysis for the proof of Theorem \ref{42thm} is a straightforward application of  the over-determined PDE machinery. Due to the size of algebraic expressions involved, the computation was performed  on  the computer algebra system \texttt{Maple}. Let us record the relevant steps of the analysis,
only for the case of Darboux super-integrable metrics. The analysis for the constant curvature metrics is similar,
and shall be omitted.

\two
Given a  Darboux super-integrable metric $\, g_+$ on a two dimensional surface $\, M$,
let $\, FM \to M$ denote the associated  oriented orthonormal frame bundle.
Let $\, \{ \, \omega^1, \, \omega^2, \, \rho \, \}$ be the canonical 1-forms on $\, FM$
which satisfy the structure equation  \eqref{42struct}.
From the general theory,
the geodesic foliation of $\, g_+$ on $\, FM$ is defined by
the rank 2 sub-bundle
\beq
\mathcal{I}^+ = \langle \; \omega^2, \, \rho \; \rangle   \subset T^*FM.   \n
\eeq
Here $\, \omega^2$ defines the canonical contact structure on $\, FM$.

\two

\noi
\emph{Proof of Theorem \ref{42thm}}.

\textbf{Step 0.}
Let $\, g_-$ be another metric on $\, M$.
When pulled back to $\, FM$,
one may write
\beq
g_- = (\eta^1)^2 + (\eta^2)^2, \n
\eeq
for a  $\, g_-$-orthonormal coframe $\, \{ \,  \eta^1, \, \eta^2 \, \}$,
where
\begin{align}\label{421eta}
\eta^1 &= a \, \omega^1 + b  \, \omega^2,   \\
\eta^2 &=  \qquad \quad    c  \, \omega^2, \n
\end{align}
for coefficients  $\, a, \, b, \, c$.
Let $\, \psi$ denote  the connection 1-form for $\, \{ \,  \eta^1, \, \eta^2 \, \}$  so that
\begin{align} \label{421etastruct}
\ed \eta^1 &= -  \psi \w         \eta^2,    \\
\ed \eta^2 &= \; \; \; \psi \w  \eta^1, \n \\
\ed \psi         &= Q \, \eta^1 \w \eta^2, \n
\end{align}
where $\, Q$  is the curvature of the metric $\, g_-$.

Differentiating \eqref{421eta},
the structure equation
\eqref{421etastruct} implies that
\begin{align}
\ed a&= - b \rho \quad      \; \; \,         + a_1 \omega^1 + a_2 \omega^2, \n \\
\ed b&=  a\rho - c \psi  + b_{1} \omega^1 + b_{2} \omega^2, \n \\
\ed c&= \; \; \; b \psi   \quad \;   \,        + c_{1} \omega^1 + c_{2} \omega^2, \n
\end{align}
where $\, a_2 = b_1$.
Here $\, a_j, \, b_j, \, c_j$'s denote the covariant derivatives of $\, a, \, b, \, c$.
The connection form $\, \psi$ is given by
\beq
\psi = \frac{1}{a} (  c \rho + t \omega^1 - c_1 \omega^2), \n
\eeq
for a variable   $\, t$.
Differentiating this equation,
one may write
\beq
\ed t = 2 c_1 \rho - a_1 \psi +  t_1 \omega^1 + t_2 \omega^2. \n
\eeq
Let us remark here that
we shall adopt the similar notational convention
for the covariant derivatives, i.e.,
$\, \ed a_1 \equiv a_{11}\omega^1 + a_{12} \omega^2,  \mod \; \rho, \, \psi$, etc.
The curvature $\, Q$ for the metric  $\, g_-$ is for example given by
\beq\label{421Scurvature}
 Q = \frac{-c_{11} -t_2+ K c }{a^2 c}.
 \eeq

\one
\textbf{Step 1.}
Set  the rank 2 sub-bundles
\begin{align}
\I^0 & =  \langle \; \omega^2, \,  \omega^1 \; \rangle,    \n \\
\I^- & =  \langle \;   \eta^2, \, \psi  \; \rangle,    \n \\
       &= \langle \;  \omega^2, \, c \rho + t \omega^1  \; \rangle.  \n
\end{align}
Then  the set of  three foliations defined by
$\, \{ \,  \I^{\pm}, \, \I^0 \, \}$
is  the Legendrian 3-web $\, \W_{g_{\pm}}$ on $\, FM$,
under the non-degeneracy condition  that $\, a, \, c, \, t \ne 0$, which we assume from now on.

\one
\textbf{Step 2.}
Following the analysis of Section \ref{sec41},
it is straightforward to determine the associated coframe
$\,    \theta, \, \theta^a, \, a = 1, \, 2, \, 3$, and $\, \alpha$,    for  $\, \W_{g_{\pm}}$
which fit into the structure equation  \eqref{41struct}, \eqref{41alpha}.
The torsion coefficients $\, R, \, S, \, T, \, N, \, L$
are then expressed as the rational functions in the successive derivatives of
$\, a, \, b, \, c, \, t$.

\one
\textbf{Step 3.}
The idea is that
the integrability condition  \eqref{42const} for the Legendrian 3-web $\, \W_{g_{\pm}}$
to have maximum rank three
allows one to close up the structure equations for $\, a, \, b, \, c, \,t$.

Differentiating $\, R$ and evaluating modulo $\, \alpha, \, \omega^2$,
one may solve for $\, \{ \, c_{21}, \, t_{21} \, \}$ in terms of the rest of the variables.
Differentiating $\, S$ and evaluating modulo $\, \alpha, \, \omega^2, \, \omega^1$
with these relations,
one may solve for  $\, \{ \,  b_2 \, \} $.
Successively differentiating this equation for $\, b_2$,
one may solve for $\, \{ \, b_{21}, \, b_{22}, \, c_2 \, \}$.
Successively differentiating again the equation for $\, c_2$,
one may solve for $\,\{ \,  b_{11}, \, c_{22}\, \}$.

Differentiating $\, R$ and evaluating modulo $\, \alpha, \, \omega^1, \, \rho$,
one may solve for  $\, \{ \,  t_{22} \, \}$.
This implies that $\, \ed R, \, \ed S \equiv 0, \mod \; \alpha$.

\one
\textbf{Step 4.}
At this step, the remaining undetermined second derivatives of  $\, a, \, b, \, c, \, t$ are
$\, a_{11}, \, c_{11}, \, t_{11}$.
Introduce the third order derivatives
$\, \ed a_{11} \equiv a_{111} \omega^1 + a_{112} \omega^2, \mod \; \rho$, etc,
for these variables.

From the exterior derivative identities $\, \ed (\ed (a_1)), \, \ed (\ed (c_1)), \, \ed (\ed (t_1)) = 0$,
one may solve for $\,\{ \,  a_{112}, \, c_{112}, \, t_{112} \, \}$.
From the exterior derivative identity $\, \ed (\ed (b_1))  = 0$,
one may solve for $\,\{ \,   c_{111}  \, \}$.
Differentiating $\, T$ and evaluating modulo $\, \alpha, \, \omega^2$
with  these relations,
one may solve for  $\, \{ \,  a_{111}, \, t_{111}  \, \}$.
The structure equations for  $\, a, \, b, \, c, \, t$  close  up at this step.
We remark that one already has $\, \ed  \alpha = 0$.

\one
\textbf{Step 5.}
A key integrability equation is obtained  from the identity
$\, \ed (\ed ( \frac{c}{a} a_{11} - 2 c_{11} )) = 0$,
which gives
\beq
3\,a t K_1 +   ( t b  + c a_1 - 2\,a c_1 ) K_2 = 0. \n
\eeq
Assuming the metric $\, g_+$ is not of constant curvature,
one may solve for $\, \{ \,  c_1 \, \}$.
Differentiating this equation for $\, c_1$,
one may solve for $\, \{ \, c_{11}  \, \}$.

\one
\textbf{Step 6.}
Differentiating $\, T$ and evaluating modulo $\, \alpha$,
one may solve for  $\, \{ \,  t_{11} \, \}$.
Differentiating this equation for $\, t_{11}$,
one may solve for $\, \{ \, t_1, \, t_2  \, \}$.
Note that
the set of  remaining independent variables at this step are
\linebreak
$\, \{ \, a, \, b, \, c, \, t; \, a_1, \, b_1, \,   a_{11}  \, \}$.

Taking the exterior derivative $\, \ed (\ed (a_{11}))  = 0$,
the  resulting single  integrability equation  factors into two parts.
The vanishing of the one  part  is equivalent to that
the curvature $\, Q$ of the metric $\, g_-$ is constant.
Assuming this is not the case,
the  vanishing of the other integrability equation
allows one to solve for $\, \{ \, b_1 \, \}$.
Differentiating this equation for $\, b_1$,
one may solve for $\, \{ \, a_{11}   \, \}$.

At this stage,
the remaining independent variables are
$\, \{ \, a, \, b, \, c, \, t; \, a_1   \, \}$,
and the structure equations for these variables are compatible. Since the metric $\, g_-$ is well defined on the surface $\, M$, and is  invariant under the $\, \textnormal{SO}_2$ action along the fibers of  $\,FM \to  M$,
from the general theory of differential equations it follows that the moduli space of these maximally geodesically  compatible mates is generically $\,5-1 = 4$ dimensional.

The analysis for the case when  the metric $\, g_-$  has   constant curvature follows from the similar analysis.
$\sq$

\subsubsection{Example}\label{sec422}
Lie classified the possible local symmetry Lie algebras,
and the representations thereof,
of projective vector fields  on the plane.
Bryant, Manno, \& Matveev used this
to give explicit local normal forms for
the two dimensional (pseudo) Riemannian metrics
admitting  a transitive algebra of  projective vector fields, \citep{Bret}.
In this sub-section,
we apply this result to give examples of
the maximally geodesically compatible pairs of Darboux super-integrable  metrics.
We shall closely follow \citep{Bret}.

\two
Let $\, (\, x, \, y\, )$ be  a  local coordinate of  $\, \R^2$.
Let $\, g_+$ be the  Darboux super-integrable metric
\beq\label{422g+}
g_+ =  e^{3 x}  dx^2 - 2 D_+  e^x dy^2.
\eeq
The projective connection associated with $\, g_+$,
or equivalently the equation of   un-parameterized geodesics of $\, g_+$,
is given  by the second order ODE
\beq\label{422g+ODE}
y^{''} = \frac{1}{2} y'  +D_+ e^{-2x} (y')^3.
\eeq
The local  symmetry algebra of   \eqref{422g+ODE} is generated by
\beq\label{422g+proj}
\left\langle \;
\frac{\partial}{\partial y}, \;
\frac{\partial}{\partial x}        +  y \frac{\partial}{\partial y}, \;
2y  \frac{\partial}{\partial x}   +  y^2 \frac{\partial}{\partial y} \; \right\rangle.
\eeq

One may verify by direct computation that
the projective connections invariant under \eqref{422g+proj}
are defined by the second order ODE's of the following form
\beq\label{422g-ODE}
y^{''} = \frac{1}{2} y'  + D   e^{-2x} (y')^3, \quad D \; \textnormal{is constant}.
\eeq
This projective connection is flat whenever $\, D = 0$.\footnotemark
\footnotetext{
This shows in particular that
the five parameter family of   metrics of constant curvature
which are  maximally geodesically compatible with the given metric $\, g_+$
are all  geodesically equivalent.} Since we are considering
the pair of  Darboux super-integrable  metrics, assume  $\, D \ne 0$.

The maximally geodesically compatible mates of $\, g_+$
share the   local  symmetry algebra of  projective vector fields.
It thus  suffices to find  the metrics which has \eqref{422g-ODE}
as the equation of  un-parameterized geodesics.
Consider for  example the metrics of the form
\beq
g_- = E(x, y) dx^2 + G(x, y) dy^2. \n
\eeq
Then $\, g_-$ is a mate of $\,  g_+$ with the equation of geodesics  \eqref{422g-ODE}
whenever
\begin{align}
E(x, y) & = E(x) =
\frac{  e^{x+ c_1}  } { ( - 2    D  e^{-x+c_1}  + c_2  ) ^2 },  \n \\
G(x, y) & = G(x) =
 \frac{1}{ - 2 D e^{-x+c_1} + c_2 }, \quad  \quad D \ne  D_+; \,
                     \, c_1, \, c_2 \; \mbox{are constants}. \n
\end{align}
The three (un-differentiated) Abelian relations for the Legendrian 3-web $\, \W_{g_{\pm}}$ are given by;
\begin{align}
&\Big(
 {\frac {-2 D_+{e^{-x}}{y}^{2}{p}^{2}+4 {p}^{2}{e^{x}}+{e^{x}}{y}^{2}
 -4 p{e^{x}}y}{2{p}^{2} ( D - D_+  ) }}
, \; \;
{\frac {2 D {e^{-x}}{y}^{2}{p}^{2} -4 {p}^{2}{e^{x}}-{e^{x}}{y}^{2}+4 p{e^{x}}y}
{2{p}^{2} ( D - D_+  ) }}
,\;\;
- {e^{-x}}{y}^{2} \Big)
, \n \\
&\Big(
 {\frac {-2 D_+ {e^{-x}}y{p}^{2}-2 p{e^{x}}+  {e^{x}}y}
 {2{p}^{2} ( D - D_+  ) }}
, \;\qquad \; \qquad   \quad
 {\frac {2 D {e^{-x}}y{p}^{2}+2 p{e^{x}} - {e^{x}}y  }
 {2{p}^{2} ( D - D_+  ) }}
,\; \qquad \qquad  \; \; \;  \; \;
- {e^{-x}}y \, \,      \Big)
, \n \\
&\Big(
 {\frac {-2 D_+{e^{-x}}{p}^{2}  + {e^{x}} }{2{p}^{2} ( D - D_+  ) } }
, \; \qquad  \qquad \qquad \qquad \;  \;  \quad
 {\frac { 2 D {e^{-x}}{p}^{2}  - {e^{x}}}{2{p}^{2} ( D - D_+  ) }}
,\; \qquad \; \;   \qquad \qquad       \qquad \quad\,
  {- e^{-x}} \,   \;  \,  \,   \Big)
. \n
\end{align}

%-------------------------------------------------------
%-------------------------------------------------------
\section{Legendrian $d$-webs of rank $\geq\rho_d$}\label{sec23}
We now turn our attention to the Legendrian $d$-webs for general $d\geq 3$.

\two
As the analysis in  Section \ref{sec3} will show, the essential step in establishing the rank bound for a Legendrian $d$-web is to understand the layers of compatibility equations for the linear differential system for Abelian relations.
A direct analysis of these equations, although linear, is evolved and it is difficult to draw any meaningful conclusions. 

In this section we consider the class of simple Legendrian webs defined by the second order ODE's  \eqref{23Legweb} below.
The relevant observation is that the symbol relations of compatibility equations for the Abelian relations of the general Legendrian web is presented in this class of examples in a cleaner way.  Moreover,
a generating set of first integrals for this class of Legendrian webs are written explicitly as polynomials in the adapted local coordinates. In hindsight the aforementioned compatibility equations are reflected in the higher degree algebraic relations among them, and this leads one to conclude that there exists at least $\rho_d$ linearly independent polynomial Abelian relations, Proposition \ref{23propo}.

The analysis in Section \ref{sec3}  will  show   that these account for all the Abelian relations, and the class of Legendrian $d$-webs defined by  \eqref{23Legweb}  have rank $\rho_d$.

\two
Let $\, X = J^1(\R, \, \R)$ equipped with the canonical contact structure, Example \ref{11example1}. Let $\, (x, \, y, \, p)$ be the  adapted coordinate so that the contact line bundle is generated by $\, \Co = \langle \, \ed y-p\ed x \, \rangle$.
Consider a Legendrian $\, d$-web $\, \W$ on $\, X$
defined by the following set of   second order ODE's;
\begin{align}\label{23Legweb}
& y''    = q^a, \; \;  a = 1, \, 2, \, ... \, d,     \\
& \mbox{$q^a$'s are   distinct constants}. \n
\end{align}
The corresponding dual rank 2 sub-bundles of $\, T^*X$ are
\beq\label{23bundle}
 \I^a = \langle \, \ed y-p \ed x, \; \ed p -q^a \ed x \, \rangle,
 \;   \;  a = 1, \, 2, \, ... \, d.
\eeq

\one
We claim that  \emph{$\, \W$ has at least $\rho_d$ linearly  independent Abelian relations.}\footnotemark \,
\footnotetext{
See Section \ref{sec3} for a derivation of the formula $\, \rho_d$.
}
\begin{proposition}\label{23propo}
Let $\, X = J^1(\R, \, \R)$ be the 1-jet space of scalar functions on $\, \R$ equipped with the canonical contact structure. Let $\, \W$ be the Legendrian $\, d$-web on $\, X$ defined  by the second order ODE's \eqref{23Legweb}, $\, d \geq 3$. Then $\, \W$ has at least  $ \rho_d$ linearly  independent Abelian relations, and the rank of $\, \W$ is bounded below by $\, \rho_d$.
\end{proposition}

We present a proof of the proposition in the following five steps.
In the first four steps, we  construct $\, \rho_d$ Abelian relations.
In the last step, we show that they are linearly independent.

\two
Let us indicate here a minor technical point for the analysis in this section.
An Abelian relation of a Legendrian web defined by a set of rank 2 sub-bundles $\{\I^a\}_{a=1}^d$ is by definition a $\, d$-tuple of  closed 1-forms
$\, (\, \Omega^1, \, \Omega^2, \, ... \, \Omega^d\,)$, $\, \Omega^a \in H^0(\I^a)$, such that $\sum_b\Omega^b=0$.
Fix a  reference point $\, \x_0\in X$,
and let $\, h^a$ be the  unique local anti-derivative $\, \ed h^a = \Omega^a$
in a neighborhood of $\, \x_0$ such that $\, h^a(\x_0) = 0$.
This defines an isomorphism from $\, \A(\W)$ to the space of  local  un-differentiated  Abelian relations,
which consist  of $\, d$-tuple of respective first integrals
$\, (\, h^1, \, h^2, \, ... \, h^d\,)$ such that they vanish at $\, \x_0$, and that
$\, \sum_b h^b = 0$.
We shall freely use either of these equivalent forms as convenient from now on.

\two
\noi\emph{Proof of Proposition \ref{23propo}}.
Recall $X=J^1(\R,\R)$  with the adapted coordinate $(x,y,p)$. Choose $\x_0=(0,0,0)\in X$ for a reference point. The objects of analysis below are all algebraic (analytic) and it suffices to prove the claim locally in a neighborhood of $\x_0$.

\one
\textbf{Step 1}. Universal first integrals:

Consider the second order ODE
\beq\label{23qODE}
y^{''} = q,
\eeq
where $\, q$ is an  indeterminate constant.
We wish to construct
the first integrals  for  \eqref{23qODE}
which are polynomials in $\, q$ of the form
\beq\label{23hq}
h(q) =  h_0 + h_1  \, q    + h_2 \, q^2   +  \,  ... \,h_{m-1}\, q^{m-1},
\eeq
where $\, h_j|_{\x_0} = 0$.

Set the following three first integrals for \eqref{23qODE} which are linear in $\, q$;
\begin{align}\label{23uthree}
 u^3_0   &=   (y - p \,x) + q \, \frac{x^2}{2},   \\
 u^3_1    &=   p - q\, x,  \n  \\
 u^3_2   &=    \frac{p^2}{2} - q\, y. \n
\end{align}
\noi
Note that $\,  (u^3_1)^2   =2 q \, u^3_0 + 2 u^3_2$.

Inductively define $\, u^{m+1}_j$
for $\, m \geq  3$, \,$0 \leq j \leq 2 m - 2$, by
\begin{align}
u^{m+1}_{j} &= u^3_0 \, u^m_{j}, \qquad \qquad
j = 0, \, 1, \, ... \, 2 m - 4,  \n \\
u^{m+1}_{2m-3} &= u^3_2 \, u^m_{2 m - 5},  \n \\
u^{m+1}_{2m-2} &= u^3_2 \, u^m_{2 m - 4}.  \n
\end{align}
One has the  formula
\beq\label{221uformula}
u^{m+1}_j   = ( u^3_0 )^{j_0} \, (u^3_1)^{j_1} \,  (u^3_2)^{j_2},
             \quad     0 \leq    j    \leq 2m-2,
\eeq
where  the indices $\, (j_0,\, j_1, \, j_2)$ are  uniquely determined by
\begin{align}\label{23j}
j_0+j_1+j_2&=m- 1,   \\
j_0, \, j_2 & \geq  0, \quad   j_1 =0, \, \mbox{or}\;  1,  \n  \\
j  &=  0 \cdot  j_0   + 1 \cdot  j_1 +  2 \cdot  j_2.      \n
\end{align}
Since these  are  polynomials in $\,\{ \,  u^3_0, \, u^3_1,  \, u^3_2 \, \}$,
they are all first integrals for  \eqref{23qODE}.

\two
\textbf{Step 2}. Basic properties:

Assign the weights
$$
\textnormal{weight} (x, \,y, \,p, \,q)=(-1, 0, 1, 2)
$$
respectively. We extend the weight to the polynomials in an obvious way.
Let $\, \R[\,q\,]$ be the polynomial ring in the indeterminate $\, q$.

The following lemma follows immediately from the construction of $\, u^{m+1}_j$. We shall omit the proof.
\begin{lemma}\label{23lemma1}
$\,$

a)
Each $\, u^{m+1}_j$ is  homogenous of  weight $\, j$. For a fixed $m\geq 2$,  the set of elements
$\,    \{ \, u^{m+1}_{j} \, \}_{j=0}^{2m-2}$
is linearly independent  over $\, \R$ as $\, \R[\,q\,]$-valued functions.

b)
Each $\, u^{m+1}_j$ has degree $\, m-1$ in $\, q$, and
the set of  elements
$\,    \cup_{m=2}^{\infty} \,  \{ \, u^{m+1}_{j} \, \}_{j=0}^{2m-2}$
are  linearly independent  over $\, \R$  as $\, \R[\,q\,]$-valued functions.
\end{lemma}
\noi

\one
The following refined linear independence property will be used in \textbf{Step 4.}
\begin{lemma}\label{23lemma2}
Let $\, u=u(q)$ be a finite linear combination
of the elements in
$\,    \cup_{m=2}^{\infty} \,  \{ \, u^{m+1}_{j} \, \}_{j=0}^{2m-2}$.
If   $\, u(q_0)=0$ for a value $q_0$, then  $\, u = 0$.
\end{lemma}
\Pf
Note  that
$\,
u^{m+1}_j \equiv y^{j_0} p^{j_1} (\frac{p^2}{2} - q y)^{j_2},\mod \; \, x$,
where $\,(j_0, \, j_1, \, j_2)$ is determined by \eqref{23j}.
Evaluating at $\, q=q_0$,
the highest weight term of $\, u^{m+1}_j \mod \; x$ (now forgetting the weight of $\,q$)
is $\,y^{j_0} p^{j_1+2j_2}$.
From  \eqref{23j},  since $\, j_1=0$, or $\, 1$,
the associated index  map $\, (m, \, j) \to (\, j_0, \, j =  j_1+2j_2)$
is injective.
Lemma follows from this
by applying induction on decreasing weights.
$\sq$

\two
\textbf{Step 3}. Abelian relations:

Let $\,    \mathfrak{q}^l = (  (q^1)^l, \,   (q^2)^l, \, ... \,  (q^d)^l )$
denote a vector in $\, \R^d$ for  $\, l = 0, \, 1, \,  ...  \,d-2$.
Since the constants $\, q^a$'s are distinct, by Vandermonde  identity these vectors are linearly independent and
there exists a linearly independent set of  vectors
$\, v^{\mu} = (v^{\mu}_1, \, v^{\mu}_2, \, ... \, ,   v^{\mu}_d)$,
$\, \mu = 1, \, 2, \, ... \, d-1$, which are complimentary to \,$\{ \mathfrak{q}^l \}_{l=0}^{d-2}$ such that
\beq\label{23vector}
\langle \, v^{\mu}, \,  \mathfrak{q}^l  \, \rangle
=\sum_a  \,v^{\mu}_a \, (q^a)^l = 0,  \quad
\mbox{for} \;  \; l = 0, \, 1, \, ... \, d-\mu-1.
\eeq

Take an element $\, h = u^{m+1}_j$, $\, m \leq d-1$,
and  consider the first integral  $\, h(q^a)$  for  $\, \I^a$, \eqref{23bundle}.
By definition of  $\, u^{m+1}_j$, it  is  a polynomial of  the form
$$
 h(q^a)  = h_0 + h_1  \, q^a   + h_2 \,(q^a)^2   +  \,  ... \, ,  h_{m-1}\, (q^a)^{m-1},
$$
which is of degree $\, m-1$ in $\, q^a$.
It follows that for each $\, v^{\mu}$, $\, 1 \leq \mu \leq  d-m$,
the $\, d$-tuple of first integrals
\beq\label{23vectorh}
\Big(  v^{\mu}_1 \, h(q^1), \, v^{\mu}_2 \, h(q^2), \,  ... \, ,  \, v^{\mu}_d \,  h(q^d)\Big)
\eeq
gives an Abelian relation for $\, \W$
(the sum of components vanish by  \eqref{23vector}).

\two
In summary,
\emph{\,each $\, u^{m+1}_j$, $\, 2 \leq m \leq d-1$ and  $\,0 \leq j \leq 2m-2$,
gives rise to
$\,  d-m $ Abelian relations.}

\two
\textbf{Step 4}.   Decomposition of $\, \rho_d$:

The preceding analysis suggests the following decomposition of $\, \rho_d$.
\begin{align}\label{23decomp}
 \rho_3&: \quad 1\,\cdot 3     \\
 \rho_4&: \quad 2\,\cdot 3       +  1\,\cdot 5   \n \\
 \rho_5&: \quad 3\,\cdot 3       +  2\,\cdot 5
              +  1\,\cdot   7  \n \\
 \rho_6&: \quad 4\,\cdot 3       +  3\,\cdot 5
              +  2\,\cdot   7   +  1\,\cdot  9  \n \\
& \; \; \;  ... \n \\
 \rho_d&: \quad (d-2) \cdot 3   +  (d-3) \cdot 5   + (d-4)  \cdot   7
                  \; \; \; ... \; \; \;   +  1 \cdot  (2d-3).  \n
\end{align}
The set of  universal first integrals
$\,\{ \, u^{m+1}_{j} \, \}_{j=0}^{2m-2}$ for $\, m = 2, \, 3, \, ... \,$
are identified with  the diagonal  entries
$\, 1 \cdot  3, \, 1 \cdot  5, \, 1 \cdot  7, \,1 \cdot  9, \,  ... \,$,
in the above  decomposition.
The entries on the column below each diagonal element
are identified with those generated by
the appropriate vectors $\, v^{\mu}$, \eqref{23vector}, \eqref{23vectorh}.
One may check
that this yields the correct formula
for $\, \rho_d$.

\two
\textbf{Step 5}.   Linear independence:

From \textbf{Step 3}, for fixed $d\geq 3$
the given set of $\, \rho_d$ Abelian relations are explicitly written by
\begin{align}\label{23rhod}
&\Gamma^{m; \,\mu}_j =\Big(\, v^{\mu}_1 \,  u^{m+1}_j(q^1), \;  v^{\mu}_2 \, u^{m+1}_j(q^2), \,
... \; ,  \,
v^{\mu}_d \, u^{m+1}_j(q^d) \, \Big),   \\
&\mbox{for} \; \;  2 \leq m \leq d-1, \; \, 0 \leq j \leq 2m-2, \; \,
              1 \leq \mu \leq d-m.  \n
\end{align}
Suppose for a set of coefficients $\, c_{m;\, \mu}^j$
a linear combination
$$ \sum_{m;  \, \mu,   \,  j}   c_{m;\, \mu}^j \, \Gamma^{m; \, \mu}_j =0.$$
In component-wise this is equivalent to that for each $\, a = 1, \, 2, \, ... \, d$,
\beq
\sum_{m; \, j} \big( \sum_{\mu} c_{m;\, \mu}^j v^{\mu}_a \big) \,  u^{m+1}_j(q^a) =0. \n
\eeq
By Lemma \ref{23lemma2}, this implies
$\,
\sum_{m; \, j} \big( \sum_{\mu} c_{m;\, \mu}^j v^{\mu}_a \big) \,  u^{m+1}_j =0.
$
By Lemma \ref{23lemma1}, this then implies for each $\, a,\,m,\,j$ that
$$
  \sum_{\mu} c_{m;\, \mu}^j v^{\mu}_a  =0.
$$
The linear independence of the vectors
$\, v^{\mu} = (v^{\mu}_1, \, v^{\mu}_2, \, ... \, v^{\mu}_d)$,
\eqref{23vector},
then forces $\, c_{m;\, \mu}^j=0$.
$\sq$

\two
The class of Legendrian webs \eqref{23Legweb} discussed here will play a role in the analysis for the rank bound in the next section. The linear differential equation for Abelian relations of the general Legendrian web will be examined on this  class of Legendrian webs to show that the associated symbol is nondegenerate.
To this end, it would suffice to note the following intermediate bound on the rank of polynomial Abelian relations.
\begin{lemma}\label{lem:poly}
Let $\, \W$ be the Legendrian $\, d$-web defined  by the second order ODE's \eqref{23Legweb}, $\, d \geq 3$. Then $\, \W$ has exactly $ \rho_d$ linearly  independent Abelian relations which are polynomial in the variables $\{ y, p\}$.  
\end{lemma}
\Pf
Note by construction that the $\rho_d$ Abelian relations for $\W$ in the proof of Proposition \ref{23propo} are polynomial in the variables $\{ y, p\}$. We show that there are no other such polynomial Abelian relations.

\one
\textbf{Step 1}.
%An Abelian relation of  depth  $\, \leq  (2d-4)+1$  has weight $\, \leq   (2d-4)+1$ by definition of depth and weight. Since the weight  equals  the degree as a polynomial in the variables $\, \{ p, \, y \}$, an Abelian relation of  depth  $\, \leq   (2d-4)+1$  is a polynomial in the variables $\, \{ p, \, y \}$.
The first integrals  $\, \{ \, u^3_0(q^a), \, u^3_1(q^a) \, \}$
for the foliation defined by the second order ODE $\, y^{''}=q^a$
are  functionally independent
(note $ u^3_0(q^a) \equiv y, \,  \, u^3_1(q^a)\equiv  p, \mod x$).
One may in fact linearly solve for the variables  $\, \{ \, p, \, y \, \}$
in terms of   $\, \{ \, u^3_0(q^a), \, u^3_1(q^a) \, \}$.
A  first integral  for  the given ODE
which is a polynomial in  $\, \{ \, p, \, y \, \}$
thus must be a  polynomial (with constant coefficients)  of
$\, \{ \, u^3_0(q^a), \, u^3_1(q^a) \, \}$.

\textbf{Step  2}.
From the identity $\,  (u^3_1)^2 = 2 u^3_2 + 2 q u^3_0$,
one may choose as a basis for
the vector space of such polynomial first integrals;
\beq
 u^{m+1}_j(q^a)  =  ( u^3_0 )^{j_0} \, (u^3_1)^{j_1}  \,  (u^3_2)^{j_2}(q^a) , \n
\eeq
where
\begin{align}\label{eq:indexrel}
j_0+j_1+j_2&=m- 1 \geq  1,    \\
j_0, \, j_2 & \geq  0, \quad   j_1 =0, \, \mbox{or}\;  1,  \n \\
  j  &=  0 \cdot  j_0   + 1 \cdot  j_1 +  2 \cdot  j_2 \leq 2m-2.   \n
\end{align}
Note that  each  pair   $(m, \, j )$,  $\,  m \geq 2, \, 0 \leq j \leq 2m-2 $,
uniquely determines the indices $( j_0, \, j_1, \, j_2)$.

\textbf{Step  3}.
Set $Q=\{ q^1, q^2, \, ... \, q^d\}$.
Let $R$ be the vector space of (un-differentiated) Abelian relations of the form
$$R=\left\{ \Big(  f_q(x,y,p)=\sum_{2\leq m \atop 0\leq j \leq 2m-2}v_{m,j} u^{m+1}_j(q) \Big)_{q\in Q}
\; \Big{\vert} \; \sum_{q\in Q} f_q(x,y,p)=0 \right\}.$$

Assign the \emph{depths} by
$$\tn{depth}(x,y,p)=(1,2,1),$$
and extend them to the monomials in the variables $\{x, y, p\}$ in the obvious way.
A polynomial $f(x,y,p)$ is \emph{depth-homogeneous} of order $\delta$ when
$f(\lambda x, \lambda^2y, \lambda p)=\lambda^{\delta} f(x,y,p)$ for a scaling factor $\lambda$.
From the index relations \eqref{eq:indexrel}, note that each $u^{m+1}_j(q^a)$ is depth-homogeneous of order $2j_0+j_1+2j_2.$

It is clear that $R$ admits the depth decomposition
$$R=\oplus_{\delta=1}^{\infty} R(\delta),$$
where $R(\delta)$ is the space of depth-homogeneous Abelian relations of order $\delta$.

%joanne
\textbf{Step  4}.
The space of polynomial Abelian relations generated by \eqref{23rhod} admits the following description.

Let $(\tn{v}_q)_{q\in Q}$  be the (unique) solution to the Vandermonde equation
\begin{align}
 \sum_{q\in Q} \tn{v}_q q^s&=0, \quad s=0,1,\, ... \, d-2, \n\\
 \sum_{q\in Q} \tn{v}_q q^{d-1}&=1.\n
\end{align}
Then the general solution to the system of Vandermonde type equations
$$  \sum_{q\in Q} v_q q^s =0, \quad s=0,1, \, ... \, \mu,
$$
is given by 
$$\big( v_q= \tn{v}_q B(q)\big)_{q\in Q}, \quad\tn{where}\;\begin{cases}
B=0 \quad \tn{if}\; \mu\geq d-1, \\
B\in\R[t] \;\tn{of degree}\leq d-2-\mu\quad  \tn{if}\; \mu\leq d-2.
\end{cases}$$

Let $R_0(\delta)\subset R(\delta)$ be the subspace generated by the \emph{special} (\emph{monomial}) Abelian relations of the form
\begin{equation}\label{eq:R0}
\left\{ \Big(  f_q(x,y,p)= v_q u^{m+1}_j(q) \Big)_{q\in Q}
\; \Big{\vert} \; \sum_{q\in Q} f_q(x,y,p)=0 \right\}.
\end{equation}
The previous analysis, \eqref{23rhod}, shows that the solutions of Eq.~\eqref{eq:R0} are given by
\begin{equation}\label{eq:star}\begin{cases}
v_q=0, \quad \tn{if}\; m\geq d,  \\
v_q=\tn{v}_q B(q), \; \tn{where} \,B\in\R[t]\;\tn{of degree}\leq d-1-m,  \quad \tn{if}\;m\leq d-1.  
\end{cases} \end{equation}
From this, one finds that (see Table \ref{table1} in Section \ref{sec3})
$$\dim R_0(\delta)=\begin{cases}
(\ell+1)\left(d-(\ell+1)\right)\qquad \tn{if}\; \delta=2\ell\geq 2, \\
(\ell+1)\left(d-(\ell+2)\right)\qquad \tn{if}\;\delta=2\ell+1\geq 1.
\end{cases}$$
This  implies the formula
$$\sum_{\delta=1}^{2d-3}\dim R_0(\delta)=\rho_d.$$

\textbf{Step  5}.
It is thus left to show that $R(\delta)=R_0(\delta)$ and every polynomial Abelian relation is special.
The cases $\delta=1, 2$ can be checked by direct computation. We apply the induction argument on $\delta$.

Note first the formula for partial derivatives of $\{ u^3_0, u^3_1, u^3_2\}$ with respect to $(x, y)$.
\begin{center}
\begin{tabular}{>{$}c<{$} | >{$}c<{$} | >{$}c<{$}}
&  \partial_x  & \partial_y  \\ \hline
 u^3_0 &-u^3_1&1\\
u^3_1&-q&0\\
u^3_2& 0&-q
\end{tabular}
\end{center}
It follows that the partial derivative operators $\partial_x , \partial_y$ act on $R(\delta)$ such that
$$\partial_x: R(\delta)\to R(\delta-1), \quad \partial_y: R(\delta)\to R(\delta-2).$$

Recall from the index relations \eqref{eq:indexrel} that
$\tn{depth}(u^{m+1}_j)=
\tn{depth}\left(( u^3_0 )^{j_0} \, (u^3_1)^{j_1}  \,  (u^3_2)^{j_2}\right)=2j_0+j_1+2j_2.$
We treat the even, and odd depth cases separately.

\two\noi
[Case $\delta=2\ell$]\; 
Let
$$f_q=\sum_{j_0=0}^{\ell} v_{j_0; q}( u^3_0(q) )^{j_0}  (u^3_2(q))^{\ell-j_0},\quad q\in Q$$
be the components of an element in $R(\delta)$.

Taking $\partial_x$, one gets
$$\partial_x f_q=-\sum_{j_0=1}^{\ell}  j_0 v_{j_0; q}( u^3_0(q) )^{j_0-1}u^3_1(q)(u^3_2(q))^{\ell-j_0},\quad q\in Q, 
$$
which are the components of an element in $R(\delta-1)=R_0(\delta-1).$
Here the equality $R(\delta-1)=R_0(\delta-1)$ is from the induction hypothesis.

By \eqref{eq:star}, one finds that $j_0 v_{j_0;q}=\tn{v}_q B_{j_0}(q)$ for each $j_0\geq 1$, where $B_{j_0}\in\R[t]$ is of degree $\leq  d-2-\ell.$ Set
$$f'_q:=\sum_{j_0=1}^{\ell} v_{j_0; q}( u^3_0(q) )^{j_0}  (u^3_2(q))^{\ell-j_0},\quad q\in Q.$$
Then by \eqref{eq:star} again, $f'_q$'s are the components of an element in $R_0(\delta)$.
It is clear that the remaining terms
$$f''_q:=f_q-f'_q=v_{0; q} (u^3_2(q))^{\ell},\quad q\in Q$$
are the components of an element in $R_0(\delta)$.
This implies $R(\delta)=R_0(\delta)$.

\two\noi
[Case $\delta=2\ell+1$]\; 
Let
\begin{equation}\label{eq:fqfory}
f_q=\sum_{j_0=0}^{\ell} v_{j_0; q}( u^3_0(q) )^{j_0} u^3_1(q)  (u^3_2(q))^{\ell-j_0},\quad q\in Q
\end{equation}
be the components of an element in $R(\delta)$.

Taking $\partial_x$ and applying the identity $(u^3_1)^2=2qu^3_0+2u^3_2$, one gets
\begin{align}\label{eq:delxf}
\partial_x f_q&=-\sum_{j_0=1}^{\ell}  j_0 v_{j_0; q}( u^3_0(q) )^{j_0-1}(u^3_1(q))^2(u^3_2(q))^{\ell-j_0}
-\sum_{j_0=0}^{\ell}q v_{j_0; q}( u^3_0(q) )^{j_0} (u^3_2(q))^{\ell-j_0} \n \\
&=-\sum_{j_0=1}^{\ell} 2 j_0 v_{j_0; q}( u^3_0(q) )^{j_0-1} (u^3_2(q))^{\ell-j_0+1}
-\sum_{j_0=0}^{\ell}(2j_0+1)q v_{j_0; q}( u^3_0(q) )^{j_0} (u^3_2(q))^{\ell-j_0},\quad q\in Q. 
\end{align}
Similarly as before, they are the components of an element in $R(\delta-1)=R_0(\delta-1).$

By \eqref{eq:star}, an inductive argument for $j_0$ from $\ell$ to $1$ shows that
$v_{j_0;q}=\tn{v}_q B_{j_0}(q)$ for $\ell\geq j_0\geq 1$, where $B_{j_0}\in\R[t]$ is of degree $\leq  d-2-(\ell+1).$ Set
$$f'_q:=\sum_{j_0=1}^{\ell} v_{j_0; q}( u^3_0(q) )^{j_0} u^3_1(q)  (u^3_2(q))^{\ell-j_0},\quad q\in Q.$$
Then by \eqref{eq:star} again, for $d-2-(j_0+1+(\ell-j_0))=d-2-(\ell+1)$, $f'_q$'s are the components of an element in $R_0(\delta)$.

\one
Set the remaining terms
$$f''_q:=f_q-f'_q= v_{0; q} u^3_1(q)  (u^3_2(q))^{\ell},\quad q\in Q.$$
It is left to show that they are the components of an element in $R_0(\delta)$.

Now taking $\partial_y$ of Eq.~\eqref{eq:fqfory}, one gets
\begin{align}\label{eq:delyf}
\partial_y f_q&=+\sum_{j_0=1}^{\ell}  j_0 v_{j_0; q}( u^3_0(q) )^{j_0-1} u^3_1(q) (u^3_2(q))^{\ell-j_0}
-\sum_{j_0=0}^{\ell-1}q(\ell-j_0) v_{j_0; q}( u^3_0(q) )^{j_0}u^3_1(q) (u^3_2(q))^{\ell-j_0-1},\quad q\in Q. 
\end{align}
By the induction hypothesis, they are the components of an element in $R(\delta-2)=R_0(\delta-2).$
Considering the coefficients of the monomial $u^3_1(q) (u^3_2(q))^{\ell-1}$ and the normal form for $v_{1;q}$ proved earlier, a similar argument as above shows that $v_{0;q}=\tn{v}_q B_{0}(q)$  where $B_{0}\in\R[t]$ is of degree $\leq  d-2-(\ell+1).$ $\sq$

 \section{Maximum rank of a Legendrian web}\label{sec3}
%-------------------------------------------------------------------------------
\subsection{Overview}\label{sec32}
In sub-section \ref{sec321},  we record a theorem on the moduli of solutions to a closed linear differential system with constraints, Theorem \ref{32thm}. The theorem is well known and follows from an elementary application of Frobenius theorem. We record and emphasize it here because it is the main conceptual ingredient to our proof of the rank bound for Legendrian web.

In sub-section \ref{sec322}, we explain a rough idea of how to apply Theorem \ref{32thm} to the linear differential system for Abelian relations and obtain a upper-bound on the rank.

\subsubsection{Linear differential system with constraints}\label{sec321}
Let $X$ be a finite dimensional, simply connected manifold.  On $X$ suppose there be an $n$-by-$n$ matrix valued  ($n\geq 1$)   1-form $\phi$.  For a $\R^n$-valued function $f=(f^1, f^2, \,... \, f^n)$ on $X$,  consider the following closed\footnotemark
\footnotetext{Here 'closed' means that the derivatives of the unknown function $f$ are determined as the functions of $f$ itself and do not involve any new variables.} linear differential system associated with $\phi$;
\beq\label{32LDS}
\ed f = f \phi.
\eeq
Let $\mathcal{S}$ be the $\R$-vector space of solutions to \eqref{32LDS}.

Set
\beq\label{32Phi}
\Phi:=\ed\phi+\phi\w\phi.
\eeq
Differentiating \eqref{32LDS}, one gets the compatibility equation that a solution $f$ necessarily satisfies,
\beq\label{32compat}
f\Phi=0.
\eeq

\two Suppose the differential equation \eqref{32LDS}  is compatible
in the sense that the $n$-by-$n$ matrix valued 2-form $\Phi$
vanishes \beq\label{32Phi0} \Phi=0, \eeq and the equation
\eqref{32compat} holds identically. Then the following existence and
uniqueness theorem is well known.
\begin{theorem}[Moduli of solutions to a compatible, closed linear differential system]\label{32thm0}
Let $X$ be a finite dimensional, simply connected manifold. Consider a compatible, closed linear differential system \eqref{32LDS} as described above, for which \eqref{32Phi0} holds.  Let $\x_0\in X$ be a reference point. Then for any finite value $f_0\in\R^n$ there exists a unique solution to \eqref{32LDS} satisfying the initial value condition $f(\x_0)=f_0$. One consequently has that $\,\tn{dim}(\mathcal{S})=n$.
\end{theorem}

\one For our intended proof of the rank bound, we shall make use of two variants of the above theorem of the
following kinds. 

\two
Consider now the case when the two form $\Phi$ in \eqref{32Phi} does not vanish identically, and it imposes a set of $m$ linear constraint equations (not necessarily independent)
\beq\label{32constraint}
E^i:=  \sum_a C^i_a f^a=0,\quad i=1, 2, \, ... \, m.
\eeq
Here $(C^i_a)$ is an $m$-by-$n$ matrix valued function on $X$. It is possible that the successive derivatives of \eqref{32constraint} impose a sequence of additional linear compatibility equations. A closed linear differential system \eqref{32LDS} with constraint \eqref{32constraint} is called \emph{compatible} when this does not occur, and symbolically one has that
\beq\label{32constcompat}
\ed E^i\equiv0\mod \{ E^j\}_{j=1}^m, \quad  i=1, 2, \, ... \, m.
\eeq

The following theorem is an immediate application of Theorem \ref{32thm0}.
\begin{theorem}[Moduli of solutions to a closed linear differential system with constraints]\label{32thm}\;
 Let $X$ be a finite dimensional, simply connected manifold. Consider a closed linear differential system \eqref{32LDS} as described above, but with constraint \eqref{32constraint}.

\one
a) Suppose $\tn{rank}(C^i_a )$ is constant. Then
$$\tn{dim}(\mathcal{S})\leq n-\tn{rank}(C^i_a ).$$

b) Suppose the linear differential system with constraints is compatible and \eqref{32constcompat} holds.
Let $\x_0\in X$ be a reference point. Then for any finite value $f_0\in\R^n$ which satisfies the constraint equation  \eqref{32constraint} at $\x_0$, i.e., $$\sum_a C^i_a(\x_0) f_0^a=0,\quad i=1, 2, \, ... \, m,$$ 
there exists a unique solution to \eqref{32LDS} with constraint \eqref{32constraint} 
with the initial value  $f(\x_0)=f_0$. Suppose furthermore that $\tn{rank}(C^i_a)$ is constant. 
Then
$$\tn{dim}(\mathcal{S})=n-\tn{rank}(C^i_a).$$
\end{theorem}

The proof is by an application of Frobenius theorem. We shall omit the details, and refer the reader to \cite{Griffiths1992}.

%%%%%%%%%%%%%%%%%%%%%%%%%%%%
\subsubsection{Sketch of ideas}\label{sec322}
The argument for the proof of rank bound consists of the following four steps:
\tb{Initial problem},   \tb{Prolongation},   \tb{Closing up}, and  \tb{Non-degeneracy of symbol and rank bound}.

\two
\tb{Initial problem}.
Recall $X$ is a contact three manifold with the contact line bundle $\Co\subset T^*X$.
Let $\W$ be a Legendrian $d$-web on $X$ defined by a set of $d$ rank two sub-bundles $\I^a\supset \Co, \,a=1, 2, \, ... \, d$.
An Abelian relation is by definition a section of the direct sum bundle
$\, \oplus_{b=1}^d\I^b$ which satisfies the equations
\begin{align}
E_0&: \quad\mbox{$\sum_b$}\Omega^b=0,  \label{322E0}\\
P_1&: \quad\ed\Omega^a=0  \label{322P1}.
\end{align}
Here $E_0$ is a 0-th order compatibility equation, and $P_1$ is a 1-st order linear differential equation.

Let us use $\F_w$ to denote symbolically the $w$-th order derivatives of the section of
$\, \oplus_{b=1}^d\I^b$.\footnotemark
\footnotetext{In the actual analysis we shall not use the filtration by order but by \emph{depth} for the higher order derivatives. See Section \ref{sec31} for the details.} Then $E_0$ imposes a set of linear relations on $\F_0$, and $P_1$ imposes a set of linear relations on $\F_1\mod\F_0$.

\two
\tb{Prolongation}.
Differentiating $E_0$ using $P_1$ we get a new compatibility equation denoted by $E_1$, which imposes a set of linear relations on $\F_1\mod\F_0$.

The equation $P_1$ is a first order linear partial differential equation for the sections of $\oplus_{b=1}^d\I^b$. Let $P_2$ denote the second order equation obtained by differentiating $P_1$ once by the standard method of prolongation, \cite{BCG3}. It imposes a set of linear relations on $\F_2\mod  \{\F_1,\F_0\}$.

Differentiating $E_1$ using $P_2$ we get a new compatibility equation denoted by $E_2$, which imposes a set of linear relations on $\F_2\mod \{\F_1,\F_0\}$.

\two
Continuing in this manner, one obtains a sequence of pairs $(E_w, P_{w+1})$ of higher order compatibility equations $E_{w}$ and differential equations $P_{w+1}$.  The pair $(E_w, P_{w+1})$  is the $w$-th prolongation of $(E_0, P_{1})$.

\two
\tb{Closing up}.
Counting the number of variables in $\F_w$ and the number of equations in $E_w$, one finds that there exists $w=w(d)$ such that
$$\tn{\Big(number of variables in $\F_{w(d)}$\Big)} = \tn{\Big(number of equations in $E_{w(d)}$\Big)}.$$

It is at this point that we apply $(a)$ of Theorem \ref{32thm}. If the set of equations 
$E_{w(d)}\mod \{\F_w\}_{w=0}^{w(d)-1}$ has full rank on $\F_{w(d)}$, then one may solve for $\F_{w(d)}$ in terms of $\{\F_w\}_{w=0}^{w(d)-1}$. The linear differential system for Abelian relations closes up at this order $w(d)$. By  $(a)$ of Theorem \ref{32thm}, the rank bound is obtained by computing the rank of the set of equations $E_{w(d)}$,
$$ \tn{dim}(\A(\W))\leq \sum_{w=0}^{w(d)}\Big(\tn{number of variables in $\F_w$\Big)}
-\tn{rank}\Big(\cup_{w=0}^{w(d)}  E_{w}\Big).$$

The proof of rank bound for the Legendrian webs is now reduced to checking the rank of a set of linear compatibility equations \,$\cup_{w=0}^{w(d)}  E_{w}.$

\two
\tb{Non-degeneracy of symbol and rank bound}.
Let us consider for a moment the analogous problem for the case of planar webs for comparison, \cite{Bol}. The similar strategy as described above works without much change. It is easily checked that the relevant symbol for the set of compatibility equations in this case is a standard Vandermonde matrix. Non-degeneracy of the symbol follows and one gets the desired rank bound.

For Legendrian webs on the other hand, it turns out that the relevant symbol consists of the layers of matrices that contain blocks of Vandermonde-like sub-matrices. Although these are explicit integer matrices, we are currently not able to show directly that they are non-degenerate.

\two
It is at this point that we apply $(b)$ of Theorem \ref{32thm} and prove the non-degeneracy of symbol indirectly: we consider another related system of linear partial differential equation with the following properties;

\begin{enumerate}[\qquad a)]
\item the vector space of solutions are exactly the space of \emph{polynomial} Abelian relations for the Legendrian webs described in Section \ref{sec23},

\item the relevant symbol of compatibility equations is isomorphic to that of \,$\cup_{w=0}^{w(d)}  E_{w}$.
\end{enumerate}

The non-degeneracy of the entire set of compatibility equations  \,$\cup_{w=0}^{w(d)}  E_{w}$ follows from $(b)$ of Theorem \ref{32thm} by counting the number of dependent variables for the new linear differential system, and by observing that the space of solutions has the expected dimension $\rho_d$, which also gives the rank bound for our problem.

\two The remaining analysis is divided into two parts. \tb{Initial
problem} and \tb{Prolongation} are examined in Section \ref{sec31}.
\tb{Closing up} and \tb{Non-degeneracy of symbol and rank bound} are
examined in Section \ref{sec33}.

%-------------------------------------------------------------------------------
%-------------------------------------------------------------------------------
\subsection{Structure equation}\label{sec31}
We continue the analysis from \eqref{21ODE} of Section \ref{sec21}.
%%%%%%%%%%%%%%%%%%%%%%
\subsubsection{Initial problem}\label{sec311}
Let $\, \W$ be a Legendrian  $\, d$-web on a contact three-manifold  $\, X$.
Let  $\, (x , y,  p)$ be an  adapted local coordinate of $\, X$
so   that
$\, \W$  is    defined    by    the $\, d$ second order ODE's \eqref{21ODE}.
The corresponding rank 2 sub-bundles of $\, T^*X$ are locally generated by
\beq\label{31subbundle}
\I^a=\, \langle  \, \ed y-p \ed x, \; \ed p - q^a(x,y,p)  \ed x \, \rangle,
\; \; \;  a =1, 2, \, ... \, d,
\eeq
where $\, q^a \ne q^b$ for $\, a \ne b$.

Set
\begin{align}\label{31generate}
\theta     &= \ed y-p \, \ed x,  \\
\theta^a   &= \ed p-  q^a\, \ed x, \quad a=1, \, 2, \, ... \, d,  \n
\end{align}
and let
\beq\label{31Omega}
 \Omega^a = f^a_{10}  \theta^a+f^a_{01}  \theta, \quad a=1, \, 2, \, ... \, d,
\eeq
denote a section of $\, \I^a$
for a set of variables $\, \{ \, f^a_{10},  \, f^a_{01} \, \}_{a=1}^d$.

The vector space of Abelian relations
$\A(\W)$ is by definition
the space of solutions to the linear differential system
\eqref{22defi}, \eqref{22defi2}.
In terms of $\, \{ \, f^a_{10},  \, f^a_{01} \, \}_{a=1}^d$,
these equations are written as a system of first order linear differential equations with constraint as follows.
\begin{align}
 \mbox{$\sum$} \,f^b_{10}  =0, \; \, \mbox{$\sum$}  \,  q^b  f^b_{10} &=0,
\label{311zerothcompa}  \\
 \mbox{$\sum$}    \, f^b_{01}   =0.   \qquad \qquad  & \n\\
(\ed f^a_{10}+  (f^a_{01}+q^a_p f^a_{10}) \ed x  ) \w \theta^a
 +(\ed f^a_{01}+             q^a_y f^a_{10} \ed x  )\w \theta &=0,
\quad a=1, \, 2, \, ... \, d.  \label{311zerothpd}
\end{align}
Here $\, q^a_p, \, q^a_y$  denote the partial  derivatives.

%%%%%%%%%%%%%%%%%%%%%%
\subsubsection{Prolongation}\label{sec312}
Following the standard theory of exterior differential systems, \cite[Chapter VI]{BCG3}, we shall determine the infinite prolongation of the equations \eqref{311zerothcompa}, \eqref{311zerothpd}.

\two
By Cartan's lemma, there exist  the prolongation variables
$\, \{ \,  f^a_{20}, \, f^a_{11}, \, f^a_{02}   \, \}_{a=1}^d$
such that for any solution $\, \{ \, f^a_{10},  \, f^a_{01} \, \}_{a=1}^d$ one has that
\begin{align}\label{311thst}
\ed f^a_{10} &= f^a_{20} \theta^a + f^a_{11}\, \theta
               - ( \,f^a_{01}+ \tau^a_{10}) \, \ed x,    \\
\ed f^a_{01} &= f^a_{11} \theta^a  + f^a_{02}\, \theta
               - (\quad \, \quad  \; \; \tau^a_{01}) \, \ed x, \quad   a=1, \, 2, \, ... \, d,   \n
\end{align}
where $\, \tau^a_{10}=q^a_p f^a_{10}, \,   \tau^a_{01}=q^a_y f^a_{10}$.
Note the identity
$$  \tau^a_{10}, \, \tau^a_{01} \equiv 0  \mod \, f^a_{10}.$$

By successive application of Cartan's lemma (or equivalently since the mixed partials commute),
one may inductively define the higher order prolongation variables
$\, \{ \,  f^a_{i j}  \, \}_{a=1}^d$ for $\, i, \, j  \geq 0, \, (i \, j) \ne (0 \, 0)$,
by
\begin{align}\label{311inductdefi}
\ed f^a_{ij} &\equiv f^a_{(i+1)j} \theta^a + f^a_{i(j+1)}\, \theta,  \mod \; \ed x,
\end{align}
or equivalently
\beq\label{4derivatives}
f^a_{(i+1)j} =\frac{\partial}{\partial p} f^a_{ij},  \quad
f^a_{i(j+1)} =\frac{\partial}{\partial y} f^a_{ij}.
\eeq
Set   $\, f^a_{ij}=0$  for
$\,  i  < 0$, $\,  j  < 0$, or $\,(i\, j)=(0\,0)$.

\two
In order to facilitate the computation, let us  introduce the following three filtrations on the variables $\, f^a_{ij}$'s.
\begin{definition}\label{311weight}
Let $\, \{ \, f^a_{ij} \, \}_{a=1}^d, \;    i, \, j  \geq 0, \, (i\, j) \ne (0\,0)$,
be the sequence of prolongation variables inductively defined by \eqref{311thst}, \eqref{311inductdefi}.
For a single element $\, f^a_{ij}$, assign \emph{\tb{height, weight}}, and \emph{\tb{depth}} by
\begin{align}\label{311filtration}
\textnormal{height}(f^a_{ij}) &=  i,   \\
\textnormal{weight}(f^a_{ij}) &= i+j, \n \\
\textnormal{depth}(f^a_{ij})  &= i+2j. \n
\end{align}

For a finite linear combination  of  $\,  f^a_{ij}$'s,
weight (depth) is defined as  the maximum weight (depth)
of nonzero terms.
The weight (depth) of  an  Abelian relation   \eqref{31Omega},
or more generally a first integral,
is the maximum weight (depth) of the  associated nonzero derivatives $\, f^a_{ij}$.
\end{definition}

The partial filtration associated with the following subsets turns out to be convenient for our analysis.
\begin{definition}
The \emph{\tb{fundamental chamber}} $\, \Upsilon_{ij}$ at $\, (i  \,  j)$
is the subset
\beq\label{311upsilon}
\Upsilon_{ij}  = \{ \,  f^a_{i' j'}  \, |
\; \; \textnormal{weight} (f^a_{i'j'})   \leq i +j,
\; \, \textnormal{depth} (f^a_{i'j'})   \leq i +2 j \;  \, \}.
\eeq
The subset of elements $\, \{ \, f^a_{ij} \, \}_{a=1}^d \subset \Upsilon_{ij}$ are  the \emph{\tb{dominant  coefficients}} in that they are the elements for which the equalities hold in the definition \eqref{311upsilon}.
\end{definition}

With an abuse of notation, note the following differential relations.
\begin{align}\label{311Upsilonrelation}
\frac{\partial}{\partial p} \Upsilon_{i j } &\subset \Upsilon_{(i+1)j},   \\
\frac{\partial}{\partial y} \Upsilon_{i j } &\subset \Upsilon_{i(j+1)}, \n \\
\left(\frac{\partial}{\partial x} + p \frac{\partial}{\partial y} \right)
 \Upsilon_{i j } &\subset \Upsilon_{(i+1)j}, \n
\end{align}
(these are obvious except the third one, see \eqref{notetheidentity}).
These  inclusion relations will be used implicitly
for the rest of this section.

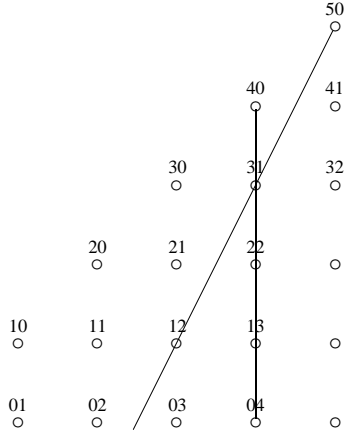
\begin{figure}[ftp]
\begin{picture}(80, 180)( -90, -10)

\multiput(160,120)(30,0){1}{$\circ$}
\multiput(190,120)(30,0){1}{$\circ$}
\multiput(190,150)(30,0){1}{$\circ$}

\multiput(130,90)(30,0){1}{$\circ$}
\multiput(160,90)(30,0){1}{$\circ$}
\multiput(190,90)(30,0){1}{$\circ$}

\multiput(100,60)(30,0){1}{$\circ$}
\multiput(130,60)(30,0){1}{$\circ$}
\multiput(160,60)(30,0){2}{$\circ$}

\multiput(70,30)(30,0){2}{$\circ$}
\multiput(130,30)(30,0){1}{$\circ$}
\multiput(160,30)(30,0){2}{$\circ$}

\multiput(70,0)(30,0){2}{$\circ$}
\multiput(130,0)(30,0){1}{$\circ$}
\multiput(160,0)(30,0){2}{$\circ$}

\multiput(69,7)(30,0){1}{\scriptsize{01}}
\multiput(69,37)(30,0){1}{\scriptsize{10}}

\multiput(99.1,67)(30,0){1}{\scriptsize{20}}
\multiput(99.1,37)(30,0){1}{\scriptsize{11}}
\multiput(99.1,7)(30,0){1}{\scriptsize{02}}

\multiput(129.3,7)(30,0){1}{\scriptsize{03}}
\multiput(129.3,37)(30,0){1}{\scriptsize{12}}
\multiput(129.3,67)(30,0){1}{\scriptsize{21}}
\multiput(129.3,97)(30,0){1}{\scriptsize{30}}

\multiput(158.9,7)(30,0){1}{\scriptsize{04}}
\multiput(158.9,37)(30,0){1}{\scriptsize{13}}
\multiput(158.9,67)(30,0){1}{\scriptsize{22}}
\multiput(158.9,97)(30,0){1}{\scriptsize{31}}
\multiput(158.9,127)(30,0){1}{\scriptsize{40}}

\multiput(188.9,127)(30,0){1}{\scriptsize{41}}
\multiput(188.9,97)(30,0){1}{\scriptsize{32}}
\multiput(188.9,157)(30,0){1}{\scriptsize{50}}

\put(162.5,4.2){\line(0,1){117}}
\put(116.1, 0){\line(1,2){75.5}}
\end{picture}
\caption{ $\, \Upsilon_{ij}$ for the case $(i  \,  j)=(3\, 1)$: upper-left chamber.}

\hspace{3.8cm}
\parbox{10cm}{
A node with $\, (i j)$ on top represents the set  $\, \{ \, f^a_{ij} \,\}_{a=1}^d$.
The weight increases as one moves to the right,
the height increases as one moves upward,
and
the depth increases as one moves to the lower-right.
The  chamber $\, \Upsilon_{ij}$
consists of the $\, f^a_{i' j'}$'s
represented by the node $\, (i'  j')$
which lie in the upper-left  chamber
bounded by the vertical weight line and the slant depth line.}
\label{311diagram}
\end{figure}

\two
The following refinement of \eqref{311inductdefi} suffices for the purpose of our analysis.
\begin{lemma}\label{311lemma}
The rough structure equation \eqref{311inductdefi}
can be refined to
\begin{align}\label{311struct}
\ed f^a_{ij} &= f^a_{(i+1)j} \theta^a + f^a_{i(j+1)}\, \theta
               - ( \, i  f^a_{(i-1)(j+1)}+ \tau^a_{ij}\,) \, \ed x,
\end{align}
where
\beq\label{31tau}
\tau^a_{ij} \equiv 0, \;  \mod \; \Upsilon_{ij}.
\eeq
\end{lemma}
\Pf
We apply induction argument on increasing depth, and weight. 

Differentiating  \eqref{311struct} by the Lie derivative
$\, L_{\frac{\partial}{\partial p}}$,
one gets
\begin{align}
\ed f^a_{(i+1) j} &= f^a_{(i+2)j} \theta^a + f^a_{(i+1)(j+1)}\, \theta
               - (i \, f^a_{i (j+1)}+ L_{\frac{\partial}{\partial p}} \tau^a_{ij}  ) \, \ed x  \n \\
               &\quad  -(  f^a_{(i+1)j} q^a_p + f^a_{i(j+1)}  )\,\ed x. \n
\end{align}
By the  induction hypothesis and \eqref{311Upsilonrelation},
$\, \{ \, f^a_{(i+1)j}, \, L_{\frac{\partial}{\partial p}} \tau^a_{ij} \, \}
\subset \Upsilon_{(i+1)j}$.
Hence
$\, \tau^a_{(i+1)j} \equiv 0,
\mod  \; \Upsilon_{(i+1)j}$.

Differentiating  \eqref{311struct} by the Lie derivative
$\, L_{\frac{\partial}{\partial y}}$,
one gets
\begin{align}
\ed f^a_{i (j+1)} &= f^a_{(i+1)(j+1)} \theta^a + f^a_{i (j+2)}\, \theta
               - (i \, f^a_{(i-1) (j+2)}+ L_{\frac{\partial}{\partial y}} \tau^a_{ij}  ) \, \ed x  \n \\
               &\quad -f^a_{(i+1)j} q^a_y  \,\ed x. \n
\end{align}
By the induction hypothesis,  one gets
$\, \tau^a_{i(j+1)} \equiv 0,  \mod  \; \Upsilon_{i(j+1)}$.

The claim \eqref{31tau} now follows from the initial   condition  \eqref{311thst},
$\, \tau^a_{10}, \, \tau^a_{01}  \equiv 0, \mod \;  \Upsilon_{10}$.
$\sq$

 %-------------------------------------------------------------------
\subsubsection{Compatibility equations}\label{sec3122}
In this sub-section, we   apply  Lemma \ref{311lemma} repeatedly
to compute the sequence of higher order compatibility equations from \eqref{311zerothcompa}.

\two
There are two preliminaries. Firstly, let us  introduce a simplified notation
\beq
\langle \mathfrak{q}^i, \, f_{jk} \rangle = \sum_b  \, (q^b)^i f^b_{jk}. \n
\eeq
For example, the initial equation \eqref{311zerothcompa} is written in this notation as
\begin{align}\label{312zerothabel}
\langle \mathfrak{q}^0, \, f_{10} \rangle&=0,
\; \, \langle \mathfrak{q}^1, \, f_{10} \rangle=0,  \\
\langle \mathfrak{q}^0, \, f_{01} \rangle&=0.  \n
\end{align}

Secondly,
set  an $\,  \infty$-by-$\infty$ integer matrix $\,  (c^I_J)$
for $\, I \geq0$ as follows.
\begin{align}\label{312cdefi}
c^I_0&=1,    \\
c^I_J&=0,  \qquad\qquad\qquad\qquad\qquad \mbox{when} \; \; J < 0,
        \; \mbox{or} \; \,  J >  \frac{I}{2}, \n \\
c^{I}_J&=  (I-2J+1)\,   c^{I-1}_{J-1}+c^{I-1}_J, \n \\
         &= \frac{I!}{2^J\,(I-2\, J)! \, J!}. \n
\end{align}
The first few terms are (for $\, J = 0, \, 1, \, 2, \, 3$),

\one

\begin{center}
\begin{tabular}{l*{5}{c}r}
  $c^0_J$    \vline& 1 & 0 & 0 & 0    \\  \hline
  $c^1_J$     \vline& 1 & 0 & 0 & 0    \\ \hline
  $c^2_J$   \vline& 1 & 1 & 0 & 0    \\   \hline
  $c^3_J$    \vline& 1 & 3 & 0 & 0    \\  \hline
  $c^4_J$   \vline& 1 & 6 & 3 & 0    \\
\end{tabular}
\end{center}

\two
The sequence of higher order compatibility equations obtained by differentiating the initial equations \eqref{311zerothcompa} is recorded in the below.  This is a result of straightforward computation. 

The symbol structure of these compatibility equations are key to establishing the rank bound. Note that the equation \eqref{312key} is written modulo lower depth terms. 
\begin{lemma}\label{312lemma}
The successive derivatives of  \eqref{311zerothcompa} by using \eqref{311struct} yield the following sequence of higher order compatibility equations.
\begin{align}\label{312key}
 E^I_{ij}: &\quad
    c^I_0 \, \langle\mathfrak{q}^I,   f_{ij} \rangle
+ c^I_1 \, \langle  \mathfrak{q}^{I-1}, f_{(i-2)(j+1)} \rangle
\, ... \,
+ c^I_k \, \langle  \mathfrak{q}^{I-k}, f_{(i-2k)(j+k)} \rangle
\, ... \,  \equiv 0,  \\
&\,\mod  \; \Upsilon_{(i+1)(j-1)},  \qquad \textnormal{for}\; \; 0 \leq I \leq i. \n
\end{align}
\end{lemma}
\Pf
We apply induction argument on
increasing depth, and weight.

Differentiating  \eqref{312key}
by  the Lie derivative    $\, L_{\frac{\partial }{\partial y }}$,
and noting  from \eqref{311Upsilonrelation} that
$\, L_{\frac{\partial }{\partial y }} \Upsilon_{(i+1)(j-1)}
\subset \Upsilon_{(i+1)j}$,
one gets
\begin{align}
 E^I_{i(j+1)}: &\quad
    c^I_0 \, \langle\mathfrak{q}^I,   f_{i(j+1)} \rangle
  + c^I_1 \, \langle  \mathfrak{q}^{I-1}, f_{(i-2)(j+2)} \rangle
\, ... \,
+ c^I_k \, \langle  \mathfrak{q}^{I-k}, f_{(i-2k)(j+k+1)} \rangle
\, ... \,  \equiv 0,    \n  \\
&\,\mod  \; \Upsilon_{(i+1)j},   \qquad   \textnormal{for}\; \; 0 \leq I \leq i. \n
\end{align}
It thus suffices to verify the case $\, j = 0$.

Differentiating  \eqref{312key}
by   the Lie derivative  $\,  L_{\frac{\partial }{\partial p}}$,
and noting from \eqref{311Upsilonrelation} that
$\, L_{\frac{\partial }{\partial p }} \Upsilon_{(i+1)(j-1)}
\subset \Upsilon_{(i+2)(j-1)}$,
one gets
\begin{align}
 E^I_{(i+1)j}: &\quad
   c^I_0 \, \langle\mathfrak{q}^I,   f_{(i+1)j} \rangle
  + c^I_1 \, \langle  \mathfrak{q}^{I-1}, f_{(i-1)(j+1)} \rangle
\, ... \,
+ c^I_k \, \langle  \mathfrak{q}^{I-k}, f_{(i-2k+1)(j+k)} \rangle
\, ... \,  \equiv 0, \n  \\
&\,\mod  \; \Upsilon_{(i+2)(j-1)},   \qquad  \textnormal{for}\; \; 0 \leq I \leq i. \n
\end{align}
It thus suffices to verify the case $\, I = i\,$.

Differentiate  \eqref{312key} now for the case  $\, I  = i$
by  the Lie derivative    $\, L_{\frac{\partial}{\partial x} + p \frac{\partial}{\partial y}}$.
One notes from \eqref{311Upsilonrelation} the relation
$\, L_{\frac{\partial}{\partial x} + p \frac{\partial}{\partial y}} \Upsilon_{(i+1)(j-1)}
\subset \Upsilon_{(i+2)(j-1)}$,
and this shows from \eqref{notetheidentity} that
\begin{align}
 E^{i+1}_{(i+1) j}:   &\quad
    c^i_0  \,   \big( \langle \mathfrak{q}^{i+1},   f_{(i+1)j} \rangle
+ i \, \langle \mathfrak{q}^{i},   f_{(i-1)(j+1) } \rangle   \big)
+   c^i_1  \,   \big( \langle \mathfrak{q}^{i},   f_{(i-1)(j+1)} \rangle
  +( i -2 ) \, \langle \mathfrak{q}^{i-1},   f_{(i-3)(j+2) } \rangle   \big)  \n \\
     &\quad ... \;  c^i_k  \,   \big( \langle \mathfrak{q}^{i+1-k},   f_{(i+1-2k)(j+k) } \rangle
+( i -2 k )\, \langle \mathfrak{q}^{i-k},   f_{( i -2 k -1)( j+k+1)} \rangle   \big)
\;  ... \,  \equiv 0, \n \\
 &\mod   \; \Upsilon_{(i+2)(j-1)}.   \n
 \end{align}
This also implies the defining recursive relation for the structural constants $c^I_J$.
\begin{align}
c^{i+1}_0&=c^{i}_0, \n \\
c^{i+1}_1&=i \, c^i_0 + c^{i}_1, \n \\
& \, ... \, \n \\
c^{i+1}_k&=(i-2k+2) \, c^i_{k-1} + c^{i}_k, \n \\
&  \, ... \,. \n
\end{align}
The claim \eqref{312key} follows from the initial  condition \eqref{312zerothabel}.
$\sq$

\two
Let us summarize the analysis so far.

\two\qquad\qquad
\parbox{13.5cm}{
$\bullet$
the sequence of higher order compatibility equations for Abelian relations $\{ E^I_{ij} \}$
propagates in the direction of increasing depth,

$\bullet$
for each fundamental chamber $\, \Upsilon_{ij}$,
there are $\, (i+1)$ linear compatibility equations
for the dominant coefficients  $\, \{ \, f^a_{ij} \, \}_{a=1}^d$.
}

\two
In the next section, we shall give a geometric interpretation of these results and give a proof of rank bound.

%-------------------------------------------------
\subsection{Proof of rank bound}\label{sec33}
In this section we give a proof that the linear differential system for Abelian relations closes up 
at depth $\, 2d-4$, and the rank of a Legendrian $\, d$-web is bounded by the expected value $\, \rho_d$.

The structure equation for Abelian relations summarized in 
Lemmas \ref{311lemma}, \ref{312lemma}  
show that the proof of rank bound 
can be reduced to the nondegeneracy of the symbol of the associated linear compatibility equations
restricted to the sequence of dominant coefficients.
In order to verify this claim, we consider the stronger truncated linear differential equation $(\star)$ in the below
for the model Legendrian webs considered in Section \ref{sec23}. 
This equation possesses the isomorphic symbols required for our analysis and yet the space of solutions 
are by definition the polynomial Abelian relations, for which the results of algebraic analysis 
from Section \ref{sec23} can be applied.
The rank bound follows from this combined with an application of Theorem \ref{32thm}.

%%%%%%%%%%%%%%%%%%%%
%%%%%%%%%%%%%%%%%%%%
\subsubsection{Model case}\label{sec331}
In this sub-section we shall consider the class of Legendrian webs
discussed in Section \ref{sec23}, which are defined by the second
order ODE's \eqref{23Legweb}. We shall use them as a model for our
indirect analysis of the symbol of linear constraint equations
\eqref{312key} via (b) of Theorem \eqref{32thm}.

\two
Let us recall some notations.
$$\left[\begin{array}{rl}
&\tn{$q^a$'s are distinct constants}, \n\\
&\theta = \ed y-p \ed x,\, \theta^a= \ed p - q^a \ed x, \n\\
&\mathfrak{q}^I = ( (q^1)^I, \, (q^2)^I, \, ... \, (q^d)^I), \n\\
&f_{ij} = ( f_{ij}^1, \, f_{ij}^2, \, ... \, f_{ij}^d).\n
\end{array}\right.$$
Consider the following linear differential system for Abelian relations which is truncated at depth $ (2d-4)+1$.
\begin{center}
\parbox{13cm}{$$(\star)\left\{
\begin{array}{rl}\label{323reduced}
\ed f^a_{ij}&= f^a_{(i+1)j} \theta^a + f^a_{i(j+1)}\, \theta
               -  \, i  f^a_{(i-1)(j+1)} \, \ed x,      \\
f^a_{ij}&=0\; \;  \textnormal{when}\; \; \,
               (i\, j)=(0\,0),  \;  \,  i  < 0, \; j  < 0, \; \tn{or}\; \,  i+2j > (2d-4)+1, \vspace{2mm}   \\
\vspace{2mm}
&\tn{with the linear constraints} \\
 E^I_{ij}&: 
    c^I_0 \, \langle \mathfrak{q}^I,   f_{ij} \rangle
+ c^I_1 \, \langle  \mathfrak{q}^{I-1}, f_{(i-2)(j+1)} \rangle
\, ... \,
+ c^I_k \, \langle  \mathfrak{q}^{I-k}, f_{(i-2k)(j+k)} \rangle
\, ... \,   =0,   \\
&\quad  \textnormal{for}\; \; 0 \leq I \leq i.
\end{array}\right.$$}
\end{center}
Let $\mathcal{S}_d$ denote the vector space of solutions to this equation.

As will be shown in the below,  
the truncation  $f^a_{ij}=0$ for $i+2j>(2d-4)+1$ in the definition of $(\star)$ implies that 
$\mathcal{S}_d$ is a priori a subspace of \emph{polynomial} Abelian relations
for the model Legendrian web.
The upshot of Lemma \ref{311lemma} and Lemma \ref{312lemma}
is that it suffices to examine the symbol of this truncated linear differential equation 
for the proof of rank bound for the general Legendrian web.

\two
For the class of Legendrian webs under consideration, one finds that 
\begin{align}
\tau^a_{ij}&=0,\; \,\tn{in}\;\; \eqref{311struct},\n\\
E^I_{ij}&\;\tn{holds}  \;  \tn{in}\;\; \eqref{312key}\;\;\tn{without}\mod\Upsilon_{(i+1)(j-1)}.\n
\end{align}
It easily follows from this that $(\star)$ is a closed linear
differential system with constraint, which is moreover compatible.
By (b) of Theorem \ref{32thm}, the rank of the set of linear
compatibility equations is determined by the number of dependent
variables $f^a_{ij}$'s and the dimension of the moduli space of
solutions.

\two
By counting for each depth, we have the following table, see Table \ref{table1}. 
\begin{table} \begin{center}
\begin{tabular}{c | c| c}
depth   & number of variables $\, f^a_{ij}$   &  number of  compatibility equations  $\, E^I_{ij}$  \\  \hline
  1     & \;\;$d$  &    $2$    \\
  2    & $2d$  &  $1+3$    \\
  3    & $2d$  &    \hspace{5.3mm} $2+4$ \\
  4    & $3d$ &   \hspace{6mm}  $1+3+5$  \\
  5     & $3d$ &   \hspace{12.1mm}  $2+4+6$      \\
  ...     & ... &   ...  \\
       &  $k d$ &   $k^2$ \\
       &  $k d$ &    $k(k+1)$ \\
    ...     & ... &   ...  \\
 $2d-4$     &  $(d-1) d$  &  $(d-1)^2$  \\
 $2d-3$     &  $(d-1) d$  &  $(d-1) d$  \\
\end{tabular} \end{center} \caption{}  \label{table1}
\end{table}
\noi
The total sum gives
\begin{align}\label{totalsum}
&\,\textnormal{\big( number of variables \big)}  \;  -  \;
\textnormal{\big( number of compatibility equations \big)}   \\
&\qquad\qquad=
\left( d +  \sum_{k=2}^{d-1}   2 k d     \right)   -
\left( 2+  \sum_{k=2}^{d-1 }  k^2 +  k(k+1)   \right) = \rho_d. \n
 \end{align}

\one
We claim that $\mathcal{S}_d$ is exactly the space of polynomial Abelian relations constructed in Section \ref{sec23}.
\begin{lemma}\label{rankboundlemma}
For each $d\geq 3$,
$$\tn{dim}\,\mathcal{S}_d=\rho_d.$$
\end{lemma}
\Pf
\iffalse
One checks that $\, (u^3_0, \, u^3_1, \, u^3_2)$ has weight   $\, (1, \, 1, \, 2)$, and  depth   $\, (2, \, 1, \, 2)$ respectively. This implies
\begin{align}\label{depthformula}
\textnormal{weight}\, ( u^{m+1}_j )  &= \;\,  j_0 + j_1 + 2 j_2,  \\
\textnormal{depth}\,   ( u^{m+1}_j )  &= 2 j_0 + j_1 + 2 j_2  \n
\end{align}
(this agrees with the previously defined notion of depth for polynomial Abelian relations in Section \ref{sec23}). 
Since $\, j_1 = 0$, or $\, 1$,  the condition
$\, \textnormal{depth}\, ( u^{m+1}_j ) \leq   (2d-4)+1$
translates  to the condition   $\, m  \leq d$.
\fi
By definition of depth and weight, an Abelian relation of  depth  $\, \leq  (2d-4)+1$  has weight $\, \leq   (2d-4)+1$. Since the weight  equals  the degree as a polynomial in the variables $\, \{ p, \, y \}$,
an Abelian relation of  depth  $\, \leq   (2d-4)+1$  is a polynomial
in the variables $\, \{ p, \, y \}$. 
The rest follows from Lemma \ref{lem:poly}. $\sq$

\two
We have the following corollaries.
\begin{corollary}\label{cor0}
The set of constraint equations $$\{ E^I_{ij} \;|\; i,j\geq 0, (i,j)\ne(0,0), i+2j\leq(2d-4)+1, \,0\leq I\leq i\}$$ has full rank for each $d\geq 3$.
\end{corollary}
\Pf This follows from \eqref{totalsum} and (b) of Theorem \ref{32thm}. $\sq$
\begin{corollary}\label{cor1}
The set of constraint equations $$\{ E^I_{ij} \;|\; i,j\geq 0, (i,j)\ne(0,0), i+2j=\delta, \,0\leq I\leq i\}$$ has full rank for each $\delta\geq 1$.
\end{corollary}
%\begin{corollary}\label{cor2}
%The class of Legendrian webs defined by the second order ODE's \eqref{23Legweb} have rank $\rho_d$.
%\end{corollary}
%\Pf From the table above, the number of variables match the number of equations at depth $2d-3$. By Corollary \ref{cor1} and \eqref{312key}, the prolongation variables $f^a_{ij}$ of depth $\geq 2d-3$ all vanish. The linear differential system for Abelian relations reduces to $(\star)$ in this case.$\sq$

%%%%%%%%%%%%%%%%%%%%%%%%
\subsubsection{Proof of rank bound}\label{sec332}
Corollary \ref{cor1}, when combined with the structure equation \eqref{311struct} and (b) of Theorem \ref{32thm}, immediately implies the desired rank bound for the Legendrian webs. Before we present a proof, let us give a possible geometric interpretation of the depth filtration and the associated layers of compatibility equations.

\two
Recall $\Co\subset T^*X$ is the contact line bundle. Given a general Legendrian web $\W$ defined by a set of rank two sub-bundles $\I^a\subset T^*X, a=1,2,\,... \, d$, the vector space of Abelian relations $\, \A(\W)$ consists of \emph{closed} sections of the direct sum bundle $\, \oplus_{a=1}^d \I^a$. Consider the projection
\beq
\pi:\A(\W) \subset H^0(\oplus_{a=1}^d \I^a)
       \to   H^0(\oplus_{a=1}^d \I^a/\Co). \n
\eeq
By the defining non-degeneracy property of a contact 1-form, 
$\pi$ is  injective. From this it follows that in order to apply Theorem \ref{32thm} to the analysis of Abelian relations, it suffices to analyze, roughly speaking, the higher order jets of sections of $\oplus_{a=1}^d \I^a$ modulo $\Co$.

In hindsight, this is due to the following observation;\;
\emph{let $ X$ be a contact three manifold with the given contact line bundle $ \Co \subset T^*X$.
Let $\, Z = (Z_1, \, Z_2, \, ... \, Z_r)$ be a $\R^r$-valued function on  $\, X$ (here $\, r$ stands for the rank of a Legendrian web). Suppose it satisfies the equation
\beq
\ed Z \equiv 0, \mod \; Z; \; \Co. \n
\eeq
Then in fact
\beq\label{keyobservation}
\ed Z \equiv 0, \mod \; Z.\one\n
\eeq}
We leave it to the reader to verify this claim.

\two
The preceding consideration can be formulated analytically as the proposed depth filtration on the higher order jets of Abelian relations.  Recall the definition of   depth$(f^a_{ij}) = i+2j$, Definition \ref{311weight}.
Set
\beq
\mathcal{F}^{\delta}  = \, \{ \, f^a_{ij} \, | \; \textnormal{depth}(f^a_{ij}) \leq  \delta  \, \}. \n
\eeq

From the structure equation in Lemma \ref{311lemma}, one finds that differentiating \emph{modulo $\Co$} corresponds to differentiating by the vector fields
$\,\frac{\partial}{\partial p}$,  $\, \frac{\partial}{\partial x} + p \frac{\partial}{\partial y}$.
One  notes the identities
\begin{align}\label{notetheidentity}
\frac{\partial}{\partial p} f^a_{ij} &\equiv  f^a_{(i+1)j},  \\
\left( \frac{\partial}{\partial x} + p \frac{\partial}{\partial y} \right)f^a_{ij}
&\equiv -q^a f^a_{(i+1)j} - i f^a_{(i-1)(j+1)},
                                     \quad  \mod \; \Upsilon_{ij}. \n
\end{align}
The operation of differentiating modulo $\Co$ increases the depth by one, and it induces a map
$$\mathcal{F}^{\delta} \to \mathcal{F}^{\delta+1}.$$

\two
Our main claim is that (symbolically),
\beq\label{33claim}
\ed   \mathcal{F}^{2d-4}  \equiv 0, \mod \; \mathcal{F}^{2d-4},
\eeq
and the linear differential system for  Abelian relations closes up at  depth  $\, 2d-4$. The desired rank bound is an immediate consequence of this claim.
\begin{theorem}\label{rankbound}
For  $d \geq 3$, the maximum rank of a Legendrian $\, d$-web is $\,\rho_d$.
\end{theorem}
\Pf
At depth $(2d-4)+1$, from Corollary \ref{cor1}, \eqref{312key}, and Table \ref{table1}, and noting  for any $i,j$ that
$$\tn{depth}(\Upsilon_{(i+1)(j-1)})\leq(i+2j)-1,$$
the set of elements $\{ f^a_{ij}\,|\, i+2j=(2d-4)+1\}$ of depth $(2d-4)+1$ can be solved in terms of $\mathcal{F}^{2d-4}$ by the set of constraint equations  $\{ E^I_{ij}\,|\, i+2j=(2d-4)+1, 0\leq I\leq i\}$. Therefore the linear differential system for Abelian relations closes up at depth $2d-4$, and all the higher depth derivatives are expressed as the (linear) functions of derivatives of depth at most $2d-4$.

In order to apply a) of Theorem \ref{32thm}, it now suffices to show that the set of remaining  constraint equations $$\{ E^I_{ij}\,|\, i+2j\leq 2d-4, 0\leq I\leq i\}$$ has full rank.

We proceed inductively on decreasing depths. At each depth $\delta\leq 2d-4$, Corollary \ref{cor1} and \eqref{312key} show that the set of constraint equations
$$\{ E^I_{ij}\,|\, i+2j =\delta, \; 0\leq I\leq i\}$$
has full rank when restricted to the set of variables of depth $\delta$ $$\{ f^a_{ij}\,|\, i+2j=\delta\}$$   modulo the lower depth variables  $\mathcal{F}^{\delta-1}$.
$\sq$
\begin{corollary}\label{cor2}
The class of Legendrian $d$-webs defined by the second order ODE's \eqref{23Legweb} have 
the maximum rank $\rho_d$.
\end{corollary}

%-------------------------------------------------
%-------------------------------------------------
\section{Concluding    remarks}\label{sec5}
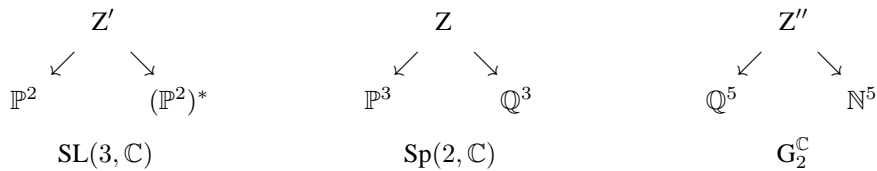
\begin{figure}[htp]
\begin{picture}(100, 75)(-80, -20)
\put(31,40){$\textnormal{Z}'$}
\put(45,25){$\searrow$}
\put(15,25){$\swarrow$}
\put(50,10){ $(\PP^2)^*$ }
\put(-4,10){ $\, \PP^2$ }

\put(18,-10){SL$(3,\C)$}

\put(161,40){$\textnormal{Z}$}
\put(175,25){$\searrow$}
\put(145,25){$\swarrow$}
\put(183,10){ $\mathbb{Q}^3$ }
\put(130,10){ $\, \PP^3$ }
%\put(185,30){$\pi_2$}
%\put(135,30){$\pi_1$}
\put(149,-10){Sp$(2,\C)$}

\put(291,40){$\textnormal{Z}''$}
\put(305,25){$\searrow$}
\put(275,25){$\swarrow$}
\put(314,10){ $\mathbb{N}^5$ }
\put(259,10){ $\,  \mathbb{Q}^5$ }

\put(290,-10){$\mbox{G}_2^{\C}$}
\end{picture}

\caption{Dualities  associated with the rank 2 simple Lie groups.}
\label{figure6}
\end{figure}

1.
The purpose of the present  paper is to propose Legendrian web as a second order  generalization of  planar web. Within the complex analytic category, the list of  dualities  associated with the rank 2 simple Lie groups in Fig.\,\ref{figure6} provides a projective geometric perspective on this generalization, \citep{Br}. Here $\, \textnormal{Z}', \, \textnormal{Z}, \, \textnormal{Z}^{''}$ are the respective incidence spaces.

From this it is a natural extension to consider the analogous problem for the webs on $\mathbb{Q}^5$
defined by a set of horizontal foliations  with respect to the $\mbox{G}_2^{\C}$ invariant rank 2 distribution (Cartan distribution).

\two 2.
(In the complex analytic category) Is there a projective geometric proof of the bound on the dimension of the space of closed holomorphic 1-forms for surfaces in $\mathbb{Q}^3$ of given degree? What is the appropriate definition of holomorphic 1-forms in this case?

The projective duality above suggests one way of defining a holomorphic 1-form on a generally singular analytic surface $\, \Sigma \subset \mathbb{Q}^3$; a meromorphic 1-form on $\Sigma$ is holomorphic when its trace on a generic small open subset of $\, \PP^3$ vanishes. Does this agree with the notion of holomorphic forms given by Henkin \& Passare in \citep{HP}? In particular, is it true that a holomorphic 1-form on a surface in this sense is necessarily closed ? (likely not)

\two 3.
Given a set of  $\, n+2\,$  points
$\, \{ \, \textnormal{p}^1, \, \textnormal{p}^2, \, ... \, \textnormal{p}^{n+2} \, \}$
in $\, \R\PP^n$ in general position,
the associated exceptional $\, n+3$ web $\, \mathcal{E}_{n+3}$
consists of the $\, n+2$ bundles of  lines  with vertices $\, \textnormal{p}^a$'s,
and
the family of  rational normal curves through   $\, \textnormal{p}^a$'s.
Damiano  in his thesis gave a characterization that, up to diffeomorphism,
$\, \mathcal{E}_{n+3}$ is  the unique non-linearizable quadrilateral web of curves
in  dimension  $\, n$, \citep{Dam}.
He also showed that
$\, \mathcal{E}_{n+3}$ is of maximum rank
based on the observation that
$\, \mathcal{E}_{n+3}$  naturally occurs on
the (smooth part of) configuration space of ordered set of  $\, n+3$ points in $\, \R\PP^1$,
which is essentially the quotient space of the Grassmannian
$\, \tn{Gr}(2, \, \R^{n+3})$
by the Cartan subgroup of  $\, \tn{SL}(n+3,\R)$.
Modulo the codimension one  subspace of  combinatorial  Abelian relations,
the remaining single non-combinatorial Abelian relation is given by
the trace of the appropriate power of the harmonic form
representing the Euler class of the canonical bundle of $\, \tn{Gr}(2, \, \R^{n+3})$.
It would be interesting if such a geometric construction exists for Legendrian webs.

\two 4.
We currently do not have any nontrivial examples of the Legendrian 4-webs of maximum rank eleven. Let us give a description of one candidate.

\two
Let  $\, (\, x^1, \, x^2, \, x^3, \, x^4\,)$ be the  coordinate of   $\, \R^4$.
Let $\, \varpi = \ed x^1 \w  \ed x^2 + \ed x^3 \w  \ed x^4$ be the standard symplectic 2-form.
Let $\, \textnormal{S}^3 \subset  \R^4$ be the unit sphere
equipped  with the $\, \textnormal{U}(2)$-invariant induced contact structure
defined by the contact 1-form
$\, \theta = -x^2 \ed x^1 + x^1 \ed x^2  -x^4 \ed x^3 + x^3 \ed x^4$.
Consider the Legendrian 4-web  $\, \W$  on a generic small  open subset of  $\, \textnormal{S}^3$
defined by the rank 2 sub-bundles
\beq
\I^a = \{ \,  \theta, \;  \ed x^a \, \}, \quad  a = 1, \, 2, \, 3, \, 4.  \n
\eeq

\one
Each of  the  four sub 3-webs of $\, \W$ has rank 0, and $\, \W$ is distinct from the examples in Section \ref{sec23}. Is $\, \W$  algebraic in any way?

%-------------------------------------------------
%-------------------------------------------------

\end{document}